\documentclass[12pt]{article}
\usepackage{amsmath}
\usepackage{amssymb}
\usepackage{cite}
\newsymbol \blackbox 1004
\newcommand{\eh}{\hfill}\newlength{\sperr}

\newenvironment{proof}{{\settowidth{\sperr}{\bf\rm
Proof}%
\par\addvspace{0.3cm}\noindent\parbox[t]{1.3\sperr}
{\bf\rm P\eh r\eh o\eh o\eh f\eh }%
}}{\nopagebreak\mbox{}
$\blackbox$\par\addvspace{0.3cm}}

\def\diag{\mathrm{diag}}
\def\a{\alpha}

\def\de{\delta}
\def\g{\gamma}

\def\k{\kappa}
\def\vk{\varkappa}

\def\s{\sigma}
\def\la{\lambda}

\def\om{\omega}

\def\d{\partial}
\def\t{\theta}

\def\Up{\Upsilon}
\def\vp{\varphi}
\def\ve{\varepsilon}
\def\wh{\widehat}
\def\wt{\widetilde}
\def\ov{\overline}

\def\BZ{{\mathbb Z}}
\def\BC{{\mathbb C}}
\def\BR{{\mathbb R}}

\def\im{{\rm Im\ }}

\newtheorem{Pa}{Paper}[section]
\newtheorem{Tm}[Pa]{{\bf Theorem}}
\newtheorem{La}[Pa]{{\bf Lemma}}
\newtheorem{Cy}[Pa]{{\bf Corollary}}
\newtheorem{Rk}[Pa]{{\bf Remark}}

\newtheorem{Ee}[Pa]{{\bf Example}}
\newtheorem{Dn}[Pa]{{\bf Definition}}
\newtheorem{Pn}[Pa]{{\bf Proposition}}
\newcommand{\CC}
{{\mathchoice {\setbox0=\hbox{$\displaystyle\rm
C$}\hbox{\hbox
to0pt{\kern0.4\wd0\vrule height0.9\ht0\hss}\box0}}
{\setbox0=\hbox{$\textstyle\rm C$}\hbox{\hbox
to0pt{\kern0.4\wd0\vrule height0.9\ht0\hss}\box0}}
{\setbox0=\hbox{$\scriptstyle\rm C$}\hbox{\hbox
to0pt{\kern0.4\wd0\vrule height0.9\ht0\hss}\box0}}
{\setbox0=\hbox{$\scriptscriptstyle\rm C$}\hbox{\hbox
to0pt{\kern0.4\wd0\vrule height0.9\ht0\hss}\box0}}}}

\title{On the GBDT version of the B\"acklund-Darboux transformation and its applications to the
linear and nonlinear equations and Weyl theory}

\author{A.L. Sakhnovich}

\date{}
\begin{document}
\maketitle

\begin{abstract}  A general theorem on the GBDT version of the B\"acklund-Darboux transformation
for systems rationally depending on the spectral parameter is treated and its applications to nonlinear
equations are given. Explicit solutions of  direct and inverse problems for Dirac-type systems, including systems with singularities, and for the system auxiliary to the $N$-wave equation are reviewed. New results on  explicit construction of the wave functions for radial
Dirac equation are obtained. 
\end{abstract}

{2010 MSC: 37K35, 34B20,  47A48, 37K10, 47A40}

{\it Keywords: B\"acklund-Darboux
transformation, Weyl function, reflection coefficient, direct problem, inverse problem,
Dirac-type system, radial Dirac equation, integrable equation}

\section{Introduction} \label{intro}
\setcounter{equation}{0}
The B\"acklund-Darboux
transformation (BDT) is  well-known in the
spectral theory and integrable nonlinear equations (see, for
instance, \cite{AvM, Ba, BC, CGHL, CC, Cie,  Da2, GeH, Gu, KSS, LRSy, LM,  MS,
Mi, RoS, ZM} and references therein).
BDT transforms initial equation or system
into another one from the same class and transforms also solutions
of the initial equation into solutions of the transformed one. 

In this paper we review results on the GBDT version of  the
B\"acklund-Darboux transformation. 
GBDT works for the matrix and scalar cases, gives  explicit expression for the iterated Darboux matrix
in terms of a transfer matrix function,  minimizes the order of the matrix that has to be inverted  
in the BDT approach. Thus, GBDT is a convenient tool to construct
wave functions and explicit solutions of the nonlinear wave equations
as well as to solve various direct and inverse problems.
GBDT and its 
applications were treated or included as important examples
in the papers \cite{FKS, FKS3, MST, SaA1, SaA2, SaA3, SaA4, SaA5, SaA6, SaA8, SaA9, SaA11, SaA12, SaA13, SaA15, SaA16, SaA17, SaA18, SZ}
(see also  \cite{GKS1, GKS2, GKS3, GKS4', GKS4, GKS6, KaS}).
Here we consider self-adjoint and skew-self-adjoint Dirac-type systems including
the singular case corresponding to soliton-positon interaction
and solve direct and inverse problems.  We solve also direct and inverse
problems for the system auxiliary to the $N$-wave equation,
construct explicit solutions of the $N$-wave equation and obtain
evolution of the corresponding Weyl function.
The results on the radial Dirac equation are new and we treat them in greater
detail.  We consider also in a detailed way a general Theorem \ref{GBDTrd}
on the GBDT for systems rationally depending on the spectral parameter
and its applications. (This result was earlier published in \cite{SaA3}.)

First, let us illustrate BDT  by the oldest and most popular example, that
is, by the Sturm-Liouville equation
\begin{equation} \label{SL}
 -\frac{d^{2}}{dx^{2}}y(x, \la)+v(x)y(x, \la)= \la y(x, \la),
\end{equation}
where
$v(x)=\ov{v}(x)$, $\, \la = \ov{\la}$. Assume that $z(x)=\ov{z}(x)$
satisfies (\ref{SL}), when $\la=c \,$, that is, $\, -z_{xx}
+v z=c z$ ($z_{xx}:=\frac{d^{2}}{dx^{2}}z$).
Then one can rewrite (\ref{SL}) in the form
\begin{equation} \label{SL1}
\Big({\cal A}^*{\cal A}+c I \Big) y(x, \la)= \la y(x, \la),
\end{equation}
where ${\cal A}$ and ${\cal A}^*$ are first order differential
expressions:
\[
{\cal A}f= \Big(\frac{d}{d x }-\frac{z_{x}}{z} \Big)f , \quad
{\cal A}^*f= -\Big(\frac{d}{d x }+ \frac{z_{x}}{z} \Big)f.
\]
Transformed equation is given by the formula
\begin{equation} \label{SL2}
\Big({\cal A}{\cal A}^*+c I \Big) y(x, \la)= \la y(x, \la).
\end{equation}
It easy to see that (\ref{SL2}) is again Sturm-Liouville equation,
but potential $v$ is transformed into $\displaystyle{ \wt
v=v-2\left(\frac{z_{x}}{z} \right)_{x}}$. Notice further
that $\Big({\cal A}{\cal A}^*+(c-\la) I \Big){\cal A}={\cal A}
\Big({\cal A}^*{\cal A}+(c-\la) I \Big)$. Hence, it follows that if
$y(x, \la)$ satisfies (\ref{SL1}), then $\wt y:= {\cal A}y$
satisfies (\ref{SL2}). Fundamental solutions of the transformed
equations can be constructed in this way. Under rather weak
conditions the spectra of operators $A^*A$ and $AA^*$ may differ
only at zero, and so under certain conditions the spectra of
Sturm-Liouville operators $L$ and $\wt L$ associated with
differential expressions $-\frac{d^2}{d x^2}+v$ and $-\frac{d^2}{d
x^2}+ \wt v$ may differ only at $c$. The B\"acklund-Darboux 
(and related commutation) methods of inserting and removing
eigenvalues of Sturm-Liouville operators
 historically go back to B\"acklund,
Darboux, and Jacobi  \cite{Ba, Da2, Jac} with decisive later contributions by
Crum, Deift, and Gesztesy \cite{Cr, D, Ge}. (See also \cite{GeT}
and a detailed account
in Appendix G in \cite{GeH}.)

One can apply BDT again to the already transformed equation
(\ref{SL2}) and so on (iterated BDT). There is also somewhat more
complicated binary BDT (see \cite{ACTV, MS}). It proves  that, if
$v$ satisfies nonlinear integrable equation, then $\wt v$ often
satisfies it too, and so BDT is used to construct solutions of the
nonlinear equations. BDT proves especially useful for the
construction of the explicit solutions, starting from the trivial initial system
(for instance, $v=0$ in (\ref{SL})).

Elementary B\"acklund-Darboux transformations for Dirac type and
more general AKNS systems one can find, for instance, in
\cite{CGHL, KoR}.  Given first order initial and  transformed
systems $u_{x}=G(x, \la)u$ and $\wt u_{x}= \wt G(x, \la)
\wt u$, their solutions are connected via so called Darboux matrix
$w$ such that $w_{x}= \wt G w-w G$. 
Here $G$, $\wt G$ and $w$ are $m \times m$ matrix functions
($m>0$).
Clearly, if $u$ satisfies
some initial system $u_{x}=G u$, then $w u$ satisfies the transformed
one $\wt u_{x}= \wt G \wt u$. Darboux matrix or gauge
transformation is of great interest in this theory (see \cite{Cie, KoR,
LRSy, MS,  SZu, TU, WF} and references therein).

We assume that  the  {\it fundamental solutions} $u$ of the considered
systems are always normalized by the condition
\begin{equation} \label{r3}
u(0,\la)=I_m,
\end{equation}
 where $I_m$  is the $m \times m$ identity matrix, 
 Then the  fundamental solution of the transformed system
is given by the formula $\wt u(x,\la)=w(x,\la)u(x,\la)w(0,\la)^{-1}$,
where $u$ is the  fundamental solution  of the initial system.

We shall  consider the GBDT version of the B\"acklund-Darboux transformation,
where the Darboux matrix is represented in the form of the transfer matrix function
$w(x,\la)=w_A(x,\la)$. Transfer matrix function corresponding to the $S$-node (transfer
matrix function in the L. Sakhnovich form) is given by the
equality
\begin{equation} \label{0.2}
w_A(\la)=  I-\Pi_2^*S^{-1}(A_1-\la I)^{-1}\Pi_1,
\end{equation}
where the matrices $A_1, \, A_2, \, S, \, \Pi_1$ and $\Pi_2$ form an $S-node$,
that is, satisfy the matrix (operator) identity 
\begin{equation} \label{r1}
A_1S-SA_2=\Pi_1\Pi_2^*.
\end{equation}
Here $S$ and $A_k$ ($k=1,2$) are $n \times n$ matrices
and $\Pi_k$ ($k=1,2)$ are $n \times m$ matrices for some integer $n>0$.
Sometimes we use notations $A_1=A$ and $\Pi_1=\Pi$.
The function $w_A$ is the transfer matrix function for system
\[
\frac{d}{dx}z=A_1z+ \Pi_1 u, \quad y=\Pi_2^*S^{-1}z+u
\]
(see Introductions in \cite{SaL1, SaL3}). At the same time
$w_A$ is a generalization of the Liv\v{s}ic-Brodskii
characteristic matrix function. 

Matrix functions
$\Pi_k(x)$ are introduced  in terms of the generalized
eigenfunction (with the generalized matrix eigenvalue) of the
initial and dual to initial systems (see Remark \ref{GEV} and
formulas  (\ref{g3})  and (\ref{g4})  for
details). Recall that eigenfunctions are essential for the
classical B\"acklund-Darboux transformation. Matrix function 
 $S(x)$ is expressed directly via $\Pi_k(x)$, and in many cases
 one can use for this purpose the identity (\ref{r1}).
Notice that  another type of matrix identities  have been successfully used for
the construction of the explicit solutions of nonlinear equations
in \cite{Mar}.
Further developments of the Marchenko
scheme are given  in  \cite{BouM, CarlS,  Sch}.  For application
of the matrix identities to the construction of solitons see also
\cite{KaGe}.

Our approach grants additionally explicit expression
for the Darboux matrix, allows to avoid the stages of the
construction of the high order matrix solutions of the nonlinear
equation and their reduction to the required order, and minimizes
the order of the matrix that has to be inverted (matrix $S$ in our
case).

GBDT for the self-adjoint and  skew-self-adjoint Dirac-type systems,
for the system auxiliary to the $N$-wave equation, GBDT for the $N$-wave
equation itself and for nonlinear Schr\"odinger equation
are studied in the next Section "Preliminaries".
A much more general Theorem \ref{GBDTrd} on GBDT and some
applications are given in Section \ref{RD}. Radial Dirac equation is treated
in Section \ref{DirRad}. Various direct and inverse problems are solved
in Section \ref{dirinv}.

By ${\cal I}$ and ${\cal I}_k$ we denote intervals on the real axis, by
diag$\{d_1, d_2, \ldots\}$ we denote diagonal matrix with the entries $d_k$
on the diagonal,  by $[D, \xi ]$ we denote the commutator $D\xi - \xi D$,
and ${\mathrm{col}}$ means column.
By the neighbourhood of zero we  mean the neighbourhood of the form $(0, \, \ve)$
or  $[0, \, \ve)$.
As usual we denote by $\BZ$ the integers, by $\BR$ the real axis, by $\BC$ the complex plain, by $\BR_+$
the positive semi-axis, and by $\BC_+$ ($\ov  \BC_+$) the open (closed) upper semi-plane.
By  $\arg(a)$ we denote the argument of $a \in \BC$.
We always assume that $\sum_{j=k}^r d_r=0$  when $k>r$. The notation $\tau \uparrow$
means that $\tau$ is a nondecreasing function.
The  spectrum of an operator $A$ is denoted by $\s(A)$.
Matrices $j$ and $J$ have the form
\begin{equation}  \label{r14}
j = \left[
\begin{array}{cc}
I_{p} & 0 \\ 0 & -I_{p}
\end{array}
\right], \quad J= \left[\begin{array}{cc}
0&I_p\\ I_p &0\end{array}\right].
\end{equation}

\section{Preliminaries} \label{Prel}
\setcounter{equation}{0}
In this section we  apply GBDT to the self-adjoint
and skew-self-adjoint Dirac-type systems and to the system auxiliary to the $N$-wave
equation. We treat these systems  on some interval
${\cal I}$ that contains $0$.
Applications to the $N$-wave and nonlinear Schr\"odinger equations
are obtained. In this way we can show how GBDT works before
the formulation of more general results.
\subsection{Gauge transformation of the Dirac-type system}  \label{subGauge}
The  self-adjoint
 Dirac-type system has the form
\begin{equation}      \label{0x1.-4}
\frac{d}{dx}u(x, \lambda )=i\big(\lambda j+jV(x)\big)u(x, \lambda ),
\end{equation}
where   $u$ is an $m \times m$ matrix function, $m=2p$,
\begin{equation}  \label{0x1.-2}
j = \left[
\begin{array}{cc}
I_{p} & 0 \\ 0 & -I_{p}
\end{array}
\right], \hspace{1em} V= \left[\begin{array}{cc}
0&v\\v^{*}&0\end{array}\right],
\end{equation}
 and the {\it potential} $v$ is a $p \times p$ matrix function.
The skew-self-adjoint Dirac-type system has the form
\begin{equation}      \label{0x1.-1}
\frac{d}{dx}u(x, \lambda )=\big(i \lambda j+j V(x)\big)u(x,
\lambda ).
\end{equation}
Dirac-type systems are also called Dirac, Zakharov-Shabat or AKNS systems. System (\ref{0x1.-4}) can be
rewritten as:
\begin{equation}       \label{se.50}
\frac{d}{dx}u(x, \lambda )+\big(\la q_1+q_0(x)\big)u(x, \lambda
)=0, \quad q_1 \equiv -i j, \quad q_0(x)=-i j V(x).
\end{equation}
GBDT for system (\ref{se.50}) is generated by
an integer $n>0$ and  three matrices, that is, by $n \times n$  matrices $A$
and $S(0)=S(0)^*$ and by  $n \times m$ matrix  $\Pi(0)$.
 It is required that these
matrices satisfy the matrix identity
\begin{equation}       \label{se.52}
AS(0)-S(0)A^*=i\Pi(0)j\Pi(0)^*.
\end{equation}
Fix $n$ and parameter matrices $A$, $S(0)$, and $\Pi(0)$.
Then we define an $n \times m$ matrix function $\Pi(x)$ by  its  value $\Pi(0)$ at $x=0$ and by the linear
differential equation
\begin{equation}       \label{se.51}
\Pi_x(x)=A \Pi(x) q_1 + \Pi(x) q_0(x).
\end{equation}
Matrix function $S(x)$ is easily recovered from  its value $S(0)$  and from the expression for its derivative:
\begin{equation}       \label{se.53}
S_x=\Pi \Pi^*.
\end{equation}
Moreover, as $S(0)=S(0)^*$ and $S_x=S_x^*$ we have $S(x)=S(x)^*$.

Recall that we consider systems and functions on the interval ${\cal I}  (i.e., $ $x \in {\cal I}$),
and that we assume  $0\in {\cal I}$. We choose point $0$, where $\Pi(x)$ and $S(x)$ are
fixed, for convenience. (It could be any other point from ${\cal I}$.)
In (\ref{se.53}) and in some of the following formulas we omit the variables
in the notations.

From (\ref{se.51}) and (\ref{se.53}) we get
\[
(AS-SA^*)_x=A\Pi \Pi^*-\Pi \Pi^*A^*,
\]
\[(i\Pi j \Pi^*)_x=(i\Pi_x
j \Pi^*)-(i\Pi_x j \Pi^*)^*=A\Pi \Pi^*-\Pi \Pi^*A^*, \,
{\mathrm{i.e.}},
\]
\begin{equation}       \label{se.54}
(AS-SA^*)_x=(i\Pi j \Pi^*)_x.
\end{equation}
Formulas (\ref{se.52}) and (\ref{se.54}) imply identity
\begin{equation}       \label{se.55}
AS(x)-S(x)A^*=i\Pi(x)j\Pi(x)^*.
\end{equation}
That is,   matrices $A$, $S(x)$, and $\Pi(x)$
satisfy identity (\ref{r1}), where
\begin{equation}       \label{r2}
 A_1=A_2^*=A, \quad \Pi_1=\Pi, \quad \Pi_2^*=i j \Pi^*,
 \end{equation}
 and we say that  $A$, $S(x)$, and $\Pi(x)$ form an $S$-node.
 Using  (\ref{r2}) rewrite formula (\ref{0.2}) for the transfer matrix function
\begin{equation}       \label{se.56}
w_A(x, \la)=  I_m-i j \Pi(x)^*S(x)^{-1}(A-\la I_n)^{-1}\Pi(x).
\end{equation}
Thus constructed matrix function $w_A$ is a gauge transformation
of the Dirac type system. To show this we shall differentiate
$w_A$ using (\ref{se.51}), (\ref{se.53}), and (\ref{se.55}). First
calculate the derivative of $\Pi^*S^{-1}$:
\begin{equation}       \label{se.57}
\Big(\Pi^*S^{-1}\Big)_x=i j \Pi^*A^*S^{-1}+i V j
\Pi^*S^{-1}-\Pi^*S^{-1}\Pi \Pi^*S^{-1}.
\end{equation}
Taking into account that (\ref{se.55}) yields $A^*S^{-1}=S^{-1}A-i
S^{-1}\Pi\Pi^*S^{-1}$, we rewrite (\ref{se.57}) as
\begin{equation}       \label{se.58}
\Big(\Pi^*S^{-1}\Big)_x=i j \Pi^*S^{-1}A+\Big(i V j+j
\Pi^*S^{-1}\Pi j-\Pi^*S^{-1}\Pi\Big) \Pi^*S^{-1}.
\end{equation}
Using (\ref{se.51}), (\ref{se.56}), and
(\ref{se.58}) one gets
\begin{eqnarray} \nonumber
\frac{d}{dx}w_A&=&-i j\Big(\big( i V j+j \Pi^*S^{-1}\Pi
j-\Pi^*S^{-1}\Pi\big)\Pi^*S^{-1}(A-\la I_n)^{-1}\Pi
\\
&& \label{se.59}
+ i j\Pi^*S^{-1}(A-\la I_n+\la I_n)(A-\la I_n)^{-1}\Pi\Big)
\\ && \nonumber
-i j \Pi^*S^{-1}(A-\la I_n)^{-1}\Big(-i(A-\la I_n+\la I_n) \Pi j-i
\Pi j V \Big).
\end{eqnarray}
From (\ref{se.56}) and (\ref{se.59}) it follows that
\begin{eqnarray} \nonumber
\frac{d}{dx}w_A&=&\big( i j V + \Pi^*S^{-1}\Pi -j \Pi^*S^{-1}\Pi j
\big)(w_A-I_n)+\Pi^*S^{-1}\Pi
\\ && \nonumber
+i \la j (w_A-I_n)-j \Pi^*S^{-1}\Pi j-(w_A-I_n)(i \la j +i j V) \\
\nonumber
&=&\big(i \la j+ i j V + \Pi^*S^{-1}\Pi -j \Pi^*S^{-1}\Pi j
\big)w_A-w_A(i \la j +i j V).
\end{eqnarray}
Finally, rewrite relation above as
\begin{equation}       \label{se.60}
\frac{d}{dx}w_A(x, \la)=i\big( \la j+j \wt V(x)\big)w_A(x, \la)-i
w_A(x, \la)\big( \la j+j  V(x)\big),
\end{equation}
where
\begin{equation}       \label{se.61}
\wt V= \left[
\begin{array}{cc} 0 &
\wt v \\ \wt v^* & 0
\end{array}
\right]=V+i(\Pi^*S^{-1}\Pi j-j\Pi^*S^{-1}\Pi),
\end{equation}
\begin{equation}       \label{se.62}
 \wt v=v-2i[I_p
\quad 0]\Pi^*S^{-1}\Pi \left[
\begin{array}{c} 0
 \\  I_p
\end{array}
\right].
\end{equation}
Thus the following proposition is proved.
\begin{Pn}\label{PnGTDS}
Let Dirac system (\ref{0x1.-4}) and parameter matrices $A$,
$\Pi(0)$, and $S(0)=S(0)^*$ be given, and let equality
(\ref{se.52}) hold. Then the matrix function $w_A$ defined by
(\ref{se.56}), where $\Pi(x)$ and $S(x)$ are obtained via formulas
(\ref{se.51}) and (\ref{se.53}), is a gauge transformation
(Darboux matrix) of the Dirac system and satisfies equation
(\ref{se.60}). The fundamental solution of the transformed system
$\frac{d}{dx}\wt u=i( \la j+j \wt V)\wt u$ is given by the formula
\begin{equation}       \label{r19}
\wt u(x,\la)=w_A(x, \la)u(x,\la)w_A(0, \la)^{-1},
\end{equation}
where $u$ is the
fundamental solution of the initial system (\ref{0x1.-4}).
\end{Pn}
Recall that we normalize fundamental solutions $u$  by the condition 
 (\ref{r3}).

The skew-self-adjoint Dirac-type system can be written as the first
relation in (\ref{se.50}), where $q_1 \equiv -i j$, $q_0(x)=-j
V(x)$, and we put
\begin{equation}       \label{r4}
 A_1=A_2^*=A, \quad \Pi_1=\Pi, \quad \Pi_2^*=i \Pi^*,
 \end{equation}
where only the last equality differs  from the last equality in (\ref{r2}).
 After substitution of our new $q_0$ we define $\Pi$ via
(\ref{se.51}):
\begin{equation}       \label{r5}
\Pi_x(x)=A \Pi(x) q_1 + \Pi(x) q_0(x)=-A\Pi(x)j-\Pi(x)jV(x).
\end{equation}
 We  define $S$ by  the equality
\begin{equation}       \label{r6}
S_x=\Pi_1q_1\Pi_2^*=\Pi j \Pi^*.
\end{equation}
Under
these conditions and under the matrix identity at $x=0$ condition
\begin{equation}       \label{r11}
AS(0)-S(0)A^*=i\Pi(0)\Pi(0)^* ,
\end{equation}
 the identity
 \begin{equation}       \label{r11'}
AS(x)-S(x)A^*=i\Pi(x)\Pi(x)^* 
\end{equation}
follows. Similar to the previous
Proposition \ref{PnGTDS} our next proposition can be proved.
\begin{Pn}\label{PnGTskDS}
Let skew-self-adjoint Dirac-type system (\ref{0x1.-1}) be given,
and let (\ref{r11}) hold.
Then the matrix function
\begin{equation}       \label{se.63}
w_A(x, \la)=  I_m-i  \Pi(x)^*S(x)^{-1}(A-\la I_n)^{-1}\Pi(x),
\end{equation}
where $\Pi$ and $S$ are defined by  (\ref{r5}) and (\ref{r6}), is a gauge transformation
of (\ref{0x1.-1}) and satisfies equation
\begin{equation}       \label{se.64}
\frac{d}{dx}w_A(x, \la)=\big(i \la j+j \wt V(x)\big)w_A(x, \la)-
w_A(x, \la)\big(i \la j+j  V(x)\big).
\end{equation}
Here we have
\begin{equation}       \label{se.65}
\wt V= \left[
\begin{array}{cc} 0 &
\wt v \\ \wt v^* & 0
\end{array}
\right]=V+\Pi^*S^{-1}\Pi -j\Pi^*S^{-1}\Pi j.
\end{equation}
\end{Pn}
Proposition \ref{PnGTskDS} admits generalization for the case of
the system auxiliary to nonlinear optics equation
\begin{equation}       \label{se.66}
\frac{d}{dx} u(x,\la)=\big( i \la D-[D, \xi(x) ]\big)u(x,\la), \quad
\xi^*=B\xi B, 
\end{equation}
where  $u$ is  an $m \times m$ matrix function and
\begin{equation}       \label{r7}
D={\mathrm{diag}}\, \{d_1,d_2, \ldots,d_m\}=D^*; \, \,
B={\mathrm{diag}}\, \{b_1,b_2, \ldots,b_m\}, \, \, b_k=\pm 1.
\end{equation}       
To present  (\ref{se.66}) as the first equality in (\ref{se.50}) we should put $q_1 \equiv -i D$ and 
$q_0(x)=[D,\xi(x)]$. Substitute these expressions for $q_1$ and $q_0$ into (\ref{se.51})
to define $\Pi$, and 
define $S$ by the equality $S_x=\Pi D B \Pi^*$ for the derivative of $S$.  The matrix identity 
 $AS(x)-S(x)A^*=i\Pi(x)B \Pi(x)^*$ easily follows after we assume $AS(0)-S(0)A^*=i\Pi(0)B \Pi(0)^*$.
\begin{Pn}\label{PnGTNW}
Let  system (\ref{se.66}) be given. Then the matrix function
\begin{equation}       \label{se.67}
w_A(x, \la)=  I_m-i B \Pi(x)^*S(x)^{-1}(A-\la I_n)^{-1}\Pi(x)
\end{equation}
is a gauge transformation of (\ref{se.66}) and satisfies equation
\begin{equation}       \label{se.68}
\frac{d}{dx}w_A(x, \la)=\big(i \la D-[D,\wt \xi(x)]\big)w_A(x, \la)-
w_A(x, \la)\big(i \la D-[D, \xi(x)]\big).
\end{equation}
Here we have
\begin{equation}       \label{se.69}
\wt \xi=\xi  -B\Pi^*S^{-1}\Pi, \quad \wt \xi^*=B \wt \xi B.
\end{equation}
\end{Pn}
Propositions \ref{PnGTDS}-\ref{PnGTNW} are particular subcases of
a general version of B\"acklund-Darboux transformation for systems
depending rationally on the spectral parameter $\la$ (see Section \ref{RD}).
\begin{Rk} \label{GEV}
When $A$ is a scalar (i.e., $n=1$), then by (\ref{se.51}) $A$ is
an eigenvalue and $\Pi$ is an eigenfunction of the dual to
(\ref{se.50}) system $\wh u_x=(\la q_1+q_0 )\wh u$. In a general
situation we call $A$ a generalized matrix eigenvalue and we call $\Pi$ 
a generalized eigenfunction.
\end{Rk}
\subsection{$N$-wave equation and gauge transformation}  \label{subNwav}
Consider $N$-wave (nonlinear optics) equation
\begin{equation}       \label{se.70}
[D,\xi_t(x,t)]-[\wh D,\xi_x(x,t) ]=\Big[ [D,\xi(x,t) ], \, [\wh
D,\xi(x,t) ]\Big] \quad \Big( \xi_t=\frac{\d}{\d t}\xi\Big)
\end{equation}
on the product $ {\cal I}_1\times {\cal I}_2$ of  intervals $ {\cal I}_1$ and ${\cal I}_2$,
where  
\begin{equation}       \label{r8}
\xi^*=B  \xi B, \quad (0,0)\in {\cal I}_1\times {\cal I}_2, \quad \wh
D={\mathrm{diag}}\, \{\wh d_1, \ldots , \wh d_m\}=\wh D^*, 
\end{equation}
and $D$ and $B$ are given by   (\ref{r7}).
Nonlinear integrable
equation (\ref{se.70}) is the compatibility condition of two
auxiliary linear systems \cite{TF, ZaMa} (see also \cite{AbH} for the case $N>3$):
\begin{eqnarray}       \label{se.71}
&& u_x(x,t,\la) =\big( i \la D-[D, \xi(x,t) ]\big)u(x,t, \la) , \\
&&  \label{r9} u_t(x,t,\la) 
=\big( i \la \wh D-[\wh D, \xi(x,t) ]\big)u(x,t, \la) .
\end{eqnarray}
Indeed, in view of (\ref{se.71}) and (\ref{r9}) we have 
\[u_{xt}=-[D,\xi_t]u+\big( i
\la D-[D, \xi ]\big)\big( i \la \wh D-[\wh D, \xi ]\big)u,
\]
 and
\[u_{tx}=-[\wh D, \xi_x ]u+\big( i \la \wh D-[\wh D, \xi ]\big)\big(
i \la D-[D, \xi ]\big)u. 
\]
Thus one can easily see that the
compatibility condition $u_{xt}=u_{tx}$ is equivalent to
(\ref{se.70}). 

To construct gauge transformation we fix $n>0$ and three parameter matrices,
that is, two $n \times n$
 matrices $A$ and $S(0,0)=S(0,0)^*$, and  an $n \times m$
 matrix $\Pi(0,0)$ such that
\begin{equation}       \label{se.74}
AS(0,0)-S(0,0)A^*=i\Pi(0,0)B \Pi(0,0)^*.
\end{equation}
Now, introduce  matrix functions $\Pi(x,t)$ and $S(x,t)$ by the
equations
\begin{equation}       \label{se.72}
\Pi_x=-i A \Pi D+\Pi [D, \xi], \quad \Pi_t=-i A \Pi \wh D+\Pi [\wh
D, \xi],
\end{equation}
\begin{equation}       \label{se.73}
S_x=\Pi D B \Pi^*, \quad S_t=\Pi \wh D B \Pi^*.
\end{equation}
Quite similar to $u_{xt}=u_{tx}$ one can show that according to
(\ref{se.70}) equations (\ref{se.72}) are compatible, i.e.,
$\Pi_{xt}=\Pi_{tx}$. By (\ref{se.72})  and (\ref{se.73})  it is
immediate that $S_{xt}=S_{tx}$. Proposition \ref{PnGTNW} implies
\begin{Pn}\label{PnGTNWeq}
Let  $m \times m$ continuously differentiable matrix function
$\xi$ ($\xi^*=B\xi B$) satisfy $N$-wave equation (\ref{se.70}) and
let matrix functions $\Pi$ and $S=S^*$ satisfy
(\ref{se.74})-(\ref{se.73}). Then in the points of invertibility
of $S$ the matrix function
\begin{equation}       \label{se.73'}
\wt \xi(x,t):=\xi(x,t)  -B\Pi(x,t)^*S(x,t)^{-1}\Pi(x,t).
\end{equation}
satisfies equation (\ref{se.70}) and an additional condition $\wt \xi^*=B \wt \xi B$.
\end{Pn}
\begin{proof}.
Equality $\wt \xi^*=B \wt \xi B$ is immediate. Now, let $u$ be the
fundamental solution of (\ref{se.71}) and (\ref{r9}), i.e.,  let $u$ satisfy
equations  (\ref{se.71}) and (\ref{r9}), and equality $u(0,0, \la)=I_m$.
(As $\xi$ satisfies (\ref{se.70}), the compatibility condition for systems
(\ref{se.71}) and (\ref{r9}) is fulfilled, and such a matrix function $u$
exists. See, for instance, formula (1.6) on p.168 in \cite{SaL3}.)
Put
\begin{equation}       \label{se.75}
w_A(x,t, \la) =  I_m-i B \Pi(x,t)^*S(x,t)^{-1}(A-\la
I_n)^{-1}\Pi(x,t),
\end{equation}
and calculate derivatives of $w_A$ with respect to  $x$ and $t$ using Proposition \ref{PnGTNW}
in both cases. Then, for  $\wt u(x,t,\la) =w_A(x,t,\la)  u(x,t,\la)$ we have
\begin{equation}       \label{se.76}
\wt u_x(x,t,\la) =\wt G(x,t,\la)  \wt u(x,t,\la) , \quad \wt
u_t(x,t,\la) =\wt F(x,t,\la)  \wt u(x,t,\la) ,
\end{equation}
where
\begin{equation}       \label{se.77}
 \quad \wt G=i \la D-[D, \wt \xi ],
\quad \wt F=i \la \wh D-[\wh D, \wt \xi ].
\end{equation}
From the continuous differentiability of $\xi$ follows  the
continuous differentiability of $\wt \xi$ in the points of invertibility of $S$. Hence, in view of
(\ref{se.76}) we get $\wt u_{xt}=\wt u_{tx}$ or equivalently $\wt
G_t-\wt F_x+[\wt G,\wt F]=0$. As we have already discussed here,
the last equality is in its turn equivalent to the $N$-wave
equation $[D,\wt\xi_t]-[\wh D,\wt\xi_x ]=\Big[ [D,\wt\xi ], \,
[\wh D,\wt\xi ]\Big]$.
\end{proof}
 In the case of the trivial initial
solution $\xi \equiv 0$ we obtain  explicit solutions of the $N$-wave
equation. Namely, putting $\Pi(0,0)=[f_1 \, \,  f_2 \, \, \ldots \, \, f_m]$
and using (\ref{se.72}) we recover $\Pi$:
\begin{equation}       \label{se.78}
\Pi(x,t)=\Big[\exp \big(-i (d_1 x+\wh d_1 t) A\big) f_1 \quad \exp
\big(-i (d_2 x+\wh d_2 t) A\big) f_2 \quad \ldots \Big].
\end{equation}
Next we recover $S$, and explicit formulas for solutions of the
$N$-wave equation are immediate from (\ref{se.73'}).
\subsection{Nonlinear Schr\"odinger
equation: \\ n-modulation solutions } \label{subNLS}
Auxiliary systems (\ref{se.71}) and (\ref{r9}) are a particular case
of the  linear differential first order systems
\begin{equation}       \label{se.80}
u_x(x,t,\la)=G(x,t,\la)u(x,t,\la), \quad u_t(x,t,\la)=F(x,t,\la)u(x,t,\la).
\end{equation}
The compatibility condition for systems (\ref{se.80}) is given by the equation
\begin{equation}       \label{r10}
G_t-F_x+[G,F]=0.
\end{equation}
When the first and second system in (\ref{se.80}) takes the form (\ref{se.71}) and (\ref{r9}), correspondingly,
the compatibility condition (\ref{r10}) is equivalent to $N$-wave equation (\ref{se.70}).

In this subsection we consider another example, namely, the well-known integrable  \cite{ZS} focusing nonlinear Schr\"odinger equation (fNLS)
\begin{equation}       \label{se.79}
2 v_t+i(v_{xx}+2vv^*v)=0.
\end{equation}
Here $v(x,t)$ is a $p\times p$ matrix function.
Equation (\ref{se.79}) is the compatibility condition of the auxiliary systems  (\ref{se.80}), where
\begin{equation}       \label{se.80'}
 G=i\la j+j V, \quad F=i \big(\la^2 j-i
\la j V-(V_x+j V^2)/2 \big),
\end{equation}
and $j$ and $V(x,t)$ are defined in (\ref{0x1.-2}). Consider again
domain $ {\cal I}_1\times {\cal I}_2$, where $(0,0)\in {\cal
I}_1\times {\cal I}_2$. Our next proposition is a particular case
of the results of Section \ref{RD} that can be proved similar to
Propositions \ref{PnGTDS}, \ref{PnGTNWeq}.
\begin{Pn}\label{PnNLS}
Let  $p \times p$  matrix function $v(x,t)$ satisfy equation
(\ref{se.79}) and be continuously differentiable together with
$v_x$. Let $n\times m$ matrix function $\Pi$ and $n \times n$
matrix function $S=S^*$ satisfy equations
\begin{equation}       \label{se.81}
AS(0,0)-S(0,0)A^*=i \Pi(0,0) \Pi(0,0)^*,
\end{equation}
\begin{equation}       \label{se.82}
\Pi_x=-i A \Pi j- \Pi j V, \quad \Pi_t=-i A^2 \Pi j- A \Pi j
V+\frac{i}{2}\Pi(V_x+j V^2),
\end{equation}
\begin{equation}       \label{se.83}
S_x=\Pi j \Pi^*, \quad S_t=A \Pi j \Pi^*+\Pi j \Pi^* A^*-i \Pi j V
\Pi^*,
\end{equation}
where $A$ is an $n \times n$ parameter matrix. Then the identity
\begin{equation}       \label{se.83.}
AS(x,t)-S(x,t)A^*=i\Pi(x,t)\Pi(x,t)^*
\end{equation}
holds in the domain $ {\cal I}_1\times {\cal I}_2$, and  the
matrix function
\begin{equation}       \label{se.83'}
\wt v=v+2[I_p \quad 0]\Pi^*S^{-1}\Pi \left[
\begin{array}{c} 0
 \\  I_p
\end{array}
\right].
\end{equation}
 satisfies (\ref{se.79}) in the points of invertibility of
$S$. Moreover, the transfer matrix function $w_A$ given by
(\ref{se.63}) satisfies equations $(w_A)_x=\wt G w_A-w_A G$,
$(w_A)_t=\wt F w_A-w_A F$, where $\wt G$ and $\wt F$ are obtained
by substitution of $\wt V=\left[
\begin{array}{lr} 0 & \wt v
 \\  \wt v^* & 0
\end{array}
\right]$ instead of $V$ into the right-hand sides of the first and
second, respectively, relations in (\ref{se.80'}).
\end{Pn}
In the previous subsection we constructed explicit solutions of
the $N$-wave equation, using trivial initial solution $\xi=0$. In
the same way the $N$-soliton solutions of the fNLS can also be
constructed \cite{GKS2, SaA1}, putting $v=0$ in (\ref{se.82}),
(\ref{se.83}), and (\ref{se.83'}). Now, we shall study a
slightly more complicated situation, that is, the case of the nontrivial
background .
\begin{Ee} \label{n-mod} Let
$p=1$, and $v=e^{-it}$. One easily checks that
\begin{equation}       \label{se.84}
u(x,t,\la) = \big( \exp (-i t j /2) \big)C_0(\la)  \exp \big( (x+\la
t)C_1(\la)  \big),
\end{equation}
where
\begin{equation}       \label{se.85}
C_0(\la) =\left[
\begin{array}{lr} 1 & 1
 \\ -i\big( \sqrt{1+\la^2}+\la \big) & i\big( \sqrt{1+\la^2}-\la
 \big)\end{array}
\right], \quad C_1(\la) =-i \big(\sqrt{1+\la^2}\big) j,
\end{equation}
satisfies (\ref{se.80}),  (\ref{se.80'}) for the case $v=e^{-it}$. Choose for
simplicity
\begin{equation}       \label{se.86}
A={\mathrm{diag}} \, \{a_1,a_2,\ldots, a_n\}, \quad a_k\not=\ov{
a_l} \quad (k,l \leq n),
\end{equation}
and put
\begin{equation}       \label{se.87}
\psi_k(x, t)=\big(u(x,t, \ov {a_k})f_k\big)^*, \quad f_k\in \BC^2,
\quad f_k\not= 0.
\end{equation}
It follows from (\ref{se.80}), (\ref{se.80'}),  and (\ref{se.87}) that
\begin{equation}       \label{se.88}
(\psi_k^*)_x=i \ov {a_k} j\psi_k^*+j V \psi_k^*, \quad
(\psi_k^*)_t=i \ov {a_k}^2 j\psi_k^*+ \ov {a_k} j V
\psi_k^*-\frac{i}{2} j V^2 \psi_k^*.
\end{equation}
Taking into account that $v_x=0$ (i.e., $V_x=0$), by (\ref{se.86}) and (\ref{se.88}) we
see that the matrix function
\begin{equation}       \label{se.89}
\Pi=\left[
\begin{array}{c}\psi_1 \\ \ldots \\ \psi_n \end{array}
\right]
\end{equation}
satisfies (\ref{se.82}). Thus $\Pi$ is obtained from (\ref{se.87})
and (\ref{se.89}), and then relations
\begin{equation}       \label{se.90}
S=\{ s_{kj} \}_{k,j=1}^n, \quad s_{kj}=i\psi_k \psi_j^*/(a_k -
\ov{a_j})
\end{equation}
follow from the matrix identity (\ref{se.83.}). In this way,
using (\ref{se.83'}) matrix function $\wt v$ is constructed
explicitly.

Fix now integers $\{ r_{1k}\}$, $\{ r_{2k} \}$ such that
$r_{1k}^2-r_{2k}^2=l_{k}^2$ ($1\leq k\leq n$), where $l_k$ are
integer, and put $a_k=ir_{1k}/r_{2k}$. Then, by (\ref{se.84}) and
(\ref{se.87}) the dependence of $\psi_k^*$ on $t$ can be expressed
in terms of functions $\exp (\pm it/2)$ and $\exp (\pm il_k
r_{1k}t/r_{2k}^2)$. Therefore, taking into account (\ref{se.83'}),
(\ref{se.89}), and (\ref{se.90}) we see that $\wt v$ is a
periodical in $t$ solution.

Under somewhat more restrictive than (\ref{se.86}) conditions
\begin{equation}       \label{se.86'}
A={\mathrm{diag}} \, \{a_1,a_2,\ldots, a_n\}, \quad \s(A)\in
\BC_+, \quad a_k\not= a_l \, (k\not= l),
\end{equation}
matrix function $S$ is invertible, and moreover we have $S>0$. To
show this we rewrite identity (\ref{se.83.}) in the form
\begin{equation}       \label{se.91}
S(A^*-\la I_n)^{-1}-(A-\la I_n)^{-1}S=i(A-\la
I_n)^{-1}\Pi\Pi^*(A^*-\la I_n)^{-1}.
\end{equation}
As $\s(A)\in \BC_+$, by the theorem on residues one represents $S$
as integrals of the right-hand side of (\ref{se.91}) on contours
in $\ov\BC_+$, and in the limit we get
\begin{equation}       \label{se.92}
S=\frac{1}{2\pi}\int_{-\infty}^{\infty}(A-\la
I_n)^{-1}\Pi\Pi^*(A^*-\la I_n)^{-1}d \la,
\end{equation}
that is, $S\geq 0$. To derive the strict inequality suppose $S g=0$,
$g=\{g_j \}_{j=1}^n \not= 0$. Hence we have $g^*(AS-SA^*)g=0$, and
so by (\ref{se.83.}) the equality $g^*\Pi\Pi^*g=0$ holds, i.e.,
$\Pi^*g=0$. Using again (\ref{se.83.}) one gets $SA^* g=0$. Now, by
induction equalities $S\big(A^*\big)^k g=0$, $\Pi^*
\big(A^*\big)^k g=0$ ($k\geq 0$) easily follow. As by our assumption
$g\not=0$, there is its entry $g_r\not= 0$. In view of the third
relation in (\ref{se.86'}) we obtain $e(r):=\{\de_{r,j} \}_{j=1}^n\in {\mathrm{span}}
\bigcup_{k=0}^{n-1}(A^*\big)^k g$, where $\de_{r,j}$ is the Kronecker-symbol. Therefore from  $\Pi^*
\big(A^*\big)^k g=0$ ($k\geq 0$) it follows  $\Pi^* e(r)=0$, that is,
$\psi_r=0$. The last equality contradicts (\ref{se.87}), and so
inequality $S>0$ is proved. Taking into account $S>0$
we see that the fNLS solutions given by (\ref{se.83'}) are
well-defined.
\end{Ee}

\section{GBDT for system  depending rationally on\\
spectral parameter and explicit  solutions of nonlinear equations}\label{RD} 
\setcounter{equation}{0}
\subsection{GBDT for system  depending rationally on $\la$}
In this section we consider GBDT for a general case of first order 
system  depending rationally on the spectral parameter $\la$:
\begin{equation}       \label{g1}
u_x=Gu, \quad G(x,\la)=-\Big(\sum_{k=0}^r \la^k q_k(x)+\sum_{s=1}^l \sum_{k=1}^{r_s}(\la -c_s)^{-k}
q_{sk}(x)\Big),
\end{equation}
$x \in {\cal I}$, where ${\cal I }$ is an interval such that $0 \in {\cal I }$, and the coefficients
$q_k(x)$ and $q_{sk}(x)$ are $m \times m$ locally integrable matrix functions.

As before we fix an integer $n>0$. Next, we fix five  matrices,
namely, $n \times n$ matrices $A_k$ ($k=1,2$) and $S(0)$, and $n \times m$ matrices
$\Pi_k(0)$ ($k=1,2$). It is required that these matrices form an $S$-node, that is,
the identity
\begin{equation}       \label{g2}
A_1S(0)-S(0)A_2=\Pi_1(0)\Pi_2(0)^*
\end{equation}
holds. Matrix functions $\Pi_k(x)$ are introduced via the coefficients from $G$:
\begin{eqnarray}       \label{g3}
&&\big(\Pi_1 \big)_x=\sum_{k=0}^r A_1^k \Pi_1 q_k+\sum_{s=1}^l \sum_{k=1}^{r_s}(A_1 -c_s I_n)^{-k}
\Pi_1 q_{sk},\\
      \label{g4}
&&\big(\Pi_2^* \big)_x=-\Big(\sum_{k=0}^r  q_k\Pi_2^*A_2^k+\sum_{s=1}^l \sum_{k=1}^{r_s}q_{sk}\Pi_2^*(A_2 -c_s I_n)^{-k}
\Big).
\end{eqnarray}
Compare  (\ref{g1}) with (\ref{g4}) to see that $\Pi_2^*$ can be viewed
as a generalized eigenfunction of the system $u_x=Gu$.

Matrix function $S(x)$ is introduced via $\frac{d}{dx}S$ by the equality
\begin{eqnarray}       \label{g5}
S_x&=&\sum_{k=1}^r\sum_{j=1}^k A_1^{k-j} \Pi_1 q_k\Pi_2^*A_2^{j-1}-\sum_{s=1}^l \sum_{k=1}^{r_s}
\sum_{j=1}^k
(A_1 -c_s I_n)^{j-k-1} \\ \nonumber && \times
\Pi_1 q_{sk}\Pi_2^*(A_2 -c_s I_n)^{-j}.
\end{eqnarray}
Equality   (\ref{g5}) is chosen so that the identity
$\Big(A_1S-SA_2\Big)_x=\Big(\Pi_1\Pi_2^*\Big)_x$ holds.
Hence, taking into account (\ref{g2}) we have
\begin{equation}       \label{g6}
A_1S(x)-S(x)A_2=\Pi_1(x)\Pi_2(x)^*, \quad x\in {\cal I}.
\end{equation}
By Theorem \ref{GBDTrd} below, the Darboux matrix for system (\ref{g1}) has the form (\ref{0.2}) :
\begin{equation}       \label{g7}
w_A(x, \la)=  I_m-\Pi_2(x)^*S(x)^{-1}(A_1-\la I_n)^{-1}\Pi_1(x).
\end{equation}
In other words, $w_A$ satisfies the equation
\begin{equation}       \label{g8}
\frac{d}{d x}w_A(x, \la)=\wt G(x,\la)w_A(x, \la)-w_A(x, \la)G(x,\la),
\end{equation}
where $\wt G$ has the same structure as $G$:
\begin{equation}       \label{g9}
\wt G(x,\la)=-\Big(\sum_{k=0}^r \la^k \wt q_k(x)+\sum_{s=1}^l \sum_{k=1}^{r_s}(\la -c_s)^{-k}
\wt q_{sk}(x)\Big).
\end{equation}
The transformed coefficients $\wt q_k$ and $\wt q_{sk}$ are given by the formulas
\begin{equation}       \label{g10}
\wt q_k=q_k-\sum_{j=k+1}^r\Big(q_jY_{j-k-1}-X_{j-k-1}q_j+\sum_{i=k+2}^jX_{j-i}q_j
Y_{i-k-2}\Big), 
\end{equation}
\begin{eqnarray}       
  \label{g11}
\wt q_{sk}&=&q_{sk}
\\  \nonumber &&
+\sum_{j=k}^{r_s}\Big(q_{sj}Y_{s,k-j-1}-X_{s,k-j-1}q_{sj}-\sum_{i=k}^jX_{s,i-j-1}q_{sj}
Y_{s,k-i-1}\Big),
\end{eqnarray}
where $X_k(x)$, $Y_k(x)$, $X_{sk}(x)$, and $Y_{sk}(x)$ are expressed in terms of the
matrices $A_k$ and matrix functions $S(x)$ and $\Pi_k(x)$:
\begin{eqnarray}       \label{g12}
&&X_k=\Pi_2^*S^{-1}A_1^k\Pi_1, \quad Y_k=\Pi_2^*A_2^kS^{-1}\Pi_1, \\
&&  \label{g13}
X_{sk}=\Pi_2^*S^{-1}(A_1-c_s I_n)^k\Pi_1, \quad Y_{sk}=\Pi_2^*(A_2-c_s I_n)^kS^{-1}\Pi_1.
\end{eqnarray}
\begin{Tm} \label{GBDTrd} \cite{SaA3} Let first order system   (\ref{g1}) and five
matrices $S(0)$, $A_k$, and $\Pi_k$  $(k=1,2)$ be given. Assume that
the identity  (\ref{g2}) holds. Then the transfer matrix function $w_A$ given by  
(\ref{g7}), where $S$ and $\Pi_k$ are determined by (\ref{g3})--(\ref{g5}),
 satisfies equation (\ref{g8}), where $\wt G$ is determined
 by the formulas (\ref{g9})--(\ref{g13}).
 \end{Tm} 
 The proof  of Theorem \ref{GBDTrd} for the case of one pole, that is,
 for $G(x,\la)=\sum_{k=-r}^r\la^kq_k$ is contained in \cite{SaA6}.
 The case of several poles $c_s$ can be treated
 precisely  in the same way. The following formula is essential for the proof
 and is also of independent interest:
\begin{equation}       \label{g14}
\Big(\Pi_2^*S^{-1}\Big)_x=-\Big(\sum_{k=0}^r  \wt q_k\Pi_2^*S^{-1}A_1^k+\sum_{s=1}^l \sum_{k=1}^{r_s}\wt q_{sk}\Pi_2^*S^{-1}(A_1 -c_s I_n)^{-k}
\Big).
\end{equation}
Formula  (\ref{g14}) means that  multiplying $\Pi_2^*$  by $S^{-1}$ from the right
we transform a generalized eigenfunction of system (\ref{g1})  into a generalized eigenfunction
of the transformed system $\wt u_x=\wt G \wt u$. (Compare formula (\ref{g14}) with 
formula (\ref{g4}).)
\begin{Rk}\label{Rkr}
It is immediate from (\ref{g10}) and  (\ref{g12}), respectively, that $\wt q_r=q_r$ and $X_0=Y_0$.
\end{Rk}
\begin{Rk}\label{Rk3.2}
If $\s(A_1)\cap \s(A_2)=\emptyset$ the matrix function $S(x)$ is uniquely defined
by the matrix identity (\ref{g6}).
\end{Rk}
\begin{Rk}\label{Rk3.4}
In the points of the invertibility of $S$ $($and for $\la \not\in \s(A_1)\cup\s(A_2))$
the matrix function $w_A(x,\la)$ is invertible.
Indeed, from  the realization (\ref{g7}) and formula (\ref{app4}) it follows that
\[
w_A(x,  \la)^{-1}=  I_m+\Pi_2(x)^*S(x)^{-1}(A^{\times}-\la I_n)^{-1}\Pi_1(x),
\]
where $A^{\times}=A_1-\Pi_1\Pi_2^*S^{-1}$. In view of  (\ref{g6}) it is immediate that
$A^{\times}=SA_2S^{-1}$. Hence, we get
\begin{equation}       \label{g27'}
w_A(x, \la)^{-1}=  I_m+\Pi_2(x)^*(A_{2}-\la I_n)^{-1}S(x)^{-1}\Pi_1(x).
\end{equation}
\end{Rk}
\subsection{Explicit  solutions of nonlinear equations}
One can apply Theorem \ref{GBDTrd} to construct solutions
of nonlinear integrable equations and corresponding wave functions 
similar to the way, in which it was done
in subsections \ref{subNwav} and  \ref{subNLS}.
For this purpose we use auxiliary linear systems for integrable nonlinear
equation:
\begin{eqnarray}       \label{g15}
&&u_x=Gu, \quad  u_t=Fu;
\\ && \label{g16}
G(x,t,\la)=-\sum_{k=0}^r \la^k q_k(x,t)-\sum_{s=1}^l \sum_{k=1}^{r_s}(\la -c_s)^{-k}
q_{sk}(x,t),
\\ && \label{g17}
F(x,t,\la)=-\sum_{k=0}^R \la^k Q_k(x,t)-\sum_{s=1}^L \sum_{k=1}^{R_s}(\la -C_s)^{-k}
Q_{sk}(x,t),
\end{eqnarray}
and zero curvature (compatibility condition) representation (\ref{r10}) of the
integrable nonlinear equation itself. We consider nonlinear equations in the domain
$(x,t) \in  {\cal I}_1 \times {\cal I}_2$ and assume $(0,0) \in  {\cal I}_1 \times {\cal I}_2$.
\begin{Rk}\label{Rk3.3} If $G$ and $F$ are
continuously differentiable and (\ref{r10}) holds, then according to formula
(1.6)  on p.168 in \cite{SaL3} there is the $m \times m$ solution $u$
of  (\ref{g15}) normalized by the condition $u(0,0,\la)=I_m$.
\end{Rk}
\begin{Rk}\label{u}
Usually we shall asume that $G$ and $F$ are
continuously differentiable and that  $u(x,t,\la)$ is the solution of  (\ref{g15}) normalized as 
in Remark \ref{Rk3.3}.
\end{Rk}
When we deal with two auxiliary linear systems, the $n \times n$ matrix functions $S$ and 
the $n \times m$ matrix functions
$\Pi_k$
depend  on two variables $x$ and $t$, and the matrix identity   (\ref{g2}) for parameter matrices
$A_k$, $\Pi_k(0)$, and $S(0)$
is substituted by the identity 
\begin{equation}       \label{g18}
A_1S(0,0)-S(0,0)A_2=\Pi_1(0,0)\Pi_2(0,0)^*
\end{equation}
for parameter matrices $A_k$, $\Pi_k(0,0)$, and $S(0,0)$.
Equations  (\ref{g3})--(\ref{g5}) should be completed by the
similar equations with respect to derivatives in $t$.
Then Theorem \ref{GBDTrd}  provides expessions for derivatives
$\big(w_A(x,t,\la)\big)_x$ and $\big(w_A(x,t,\la)\big)_t$,
where
\begin{equation}       \label{g27}
w_A(x, t, \la)=  I_m-\Pi_2(x,t)^*S(x,t)^{-1}(A_1-\la I_n)^{-1}\Pi_1(x,t).
\end{equation}
Hence,  equations 
\begin{equation}       \label{g35}
\wt u_x=\wt G \wt u, \quad  \wt u_t=\wt F\wt u, \quad \wt u(x,t,\la):=w_A(x,t,\la)u(x,t,\la)
\end{equation}
hold. Finally  in a way which is similar to the proof of 
(\ref{g6}),
one can show that 
\begin{equation}       \label{g18'}
A_1S(x,t)-S(x,t)A_2=\Pi_1(x,t)\Pi_2(x,t)^*, \quad (x,t)\in {\cal I}_1 \times {\cal I}_2.
\end{equation}

It follows from (\ref{g35}) that
\begin{equation}       \label{g35'}
\wt u_{xt}=(\wt G_t+\wt G \wt F )\wt u, \quad  \wt u_{tx}=(\wt F_x+\wt F \wt G)\wt u.
\end{equation}
If $G$ and $F$ are continuously differentiable, then $\wt G$ and $\wt F$ are continuously differentiable too.
Hence, $\wt u_{tx}=\wt u_{xt}$ and
formula  (\ref{g35'}) implies 
\begin{equation}       \label{r12}
\big(\wt G_t - \wt F_x+[\wt G, \wt F]\big)\wt u=0.
\end{equation}
By Remark \ref{Rk3.4} and identity  (\ref{g18'}) the matrix function
$w_A(x,t,\la)$ is invertible, and by Remark \ref{Rk3.3} $u$ is invertible.
Thus, $\wt u$ is invertible.
Therefore, it is immediate from (\ref{r12})   that
\begin{equation}       \label{g18''}
\wt G_t - \wt F_x+[\wt G, \wt F]=0.
\end{equation}
For the particular case of the $N$-wave equation formula (\ref{g18''}) was obtained
in the proof of Proposition \ref{PnGTNWeq} and was used there to show that $\wt \xi$
satisfies the $N$-wave equation. Now,
we proved the following general theorem.
\begin{Tm} \label{ZQE} Let $G$ and $F$ be continuously differentiable and satisfy (\ref{r10}).
Let the identity (\ref{g2}) hold and let the matrix functions $\Pi_k$ and $S$ be given
by the equations (\ref{g3})--(\ref{g5}) and by the analogs of  (\ref{g3})--(\ref{g5}) with respect to 
the variable $t$ instead of $x$,
namely, by the equations
\begin{eqnarray}       \label{g3'}
&&\big(\Pi_1 \big)_t=\sum_{k=0}^R A_1^k \Pi_1 Q_k+\sum_{s=1}^L \sum_{k=1}^{R_s}(A_1 -C_s I_n)^{-k}
\Pi_1 Q_{sk},\\
      \label{g4'}
&&\big(\Pi_2^* \big)_t=-\Big(\sum_{k=0}^R  Q_k\Pi_2^*A_2^k+\sum_{s=1}^L \sum_{k=1}^{R_s}Q_{sk}\Pi_2^*(A_2 -C_s I_n)^{-k}
\Big).
\end{eqnarray}
\begin{eqnarray}       \label{g5'}
S_t&=&\sum_{k=1}^R \sum_{j=1}^k A_1^{k-j} \Pi_1 Q_k\Pi_2^*A_2^{j-1}-\sum_{s=1}^L \sum_{k=1}^{R_s}
\sum_{j=1}^k
(A_1 -C_s I_n)^{j-k-1} \\ \nonumber && \times
\Pi_1 Q_{sk}\Pi_2^*(A_2 -C_s I_n)^{-j}.
\end{eqnarray}
Then in the points of the invertibility of $S$ the zero curvature equation 
(\ref{g18''}), where $\wt G$ and $\wt F$ are given by (\ref{g9})--(\ref{g13}), holds.
\end{Tm}
\begin{Ee}\label{MCF}
The main chiral field equation for $m \times m$ invertible matrix function $z$ has the form
\begin{equation}       \label{g19}
2z_{xt}(x,t)=z_{x}(x,t)z(x,t)^{-1}z_{t}(x,t)+ z_{t}(x,t)z(x,t)^{-1}z_{x}(x,t),
\end{equation}
and is equivalent \cite{Po, ZM0} to the compatibility
condition (\ref{r10}) of the auxiliary systems (\ref{g15}),
where 
\begin{eqnarray}       \label{g20}
&& G(x,t,\la)=-(\la-1)^{-1}q_{11}(x,t), \quad q_{11}=z_{x}z^{-1}; 
\\  \label{g21} &&
 F(x,t,\la)=-(\la+1)^{-1}Q_{11}(x,t), 
\quad Q_{11}=-z_{t}z^{-1}.
\end{eqnarray}

Let $z$ satisfy  (\ref{g19}). In view of   (\ref{g20}) and  (\ref{g21}) equations  (\ref{g3})--(\ref{g5})
take the form
\begin{equation}       \label{g22}
\big(\Pi_1\big)_x=(A_1-I_n)^{-1}\Pi_1 z_{x}z^{-1}, \quad \big(\Pi_2^*\big)_x=-z_{x}z^{-1}\Pi_2^*(A_2-I_n)^{-1},
\end{equation}
\begin{equation}
  \label{g23}
 S_x=-(A_1-I_n)^{-1}\Pi_1 z_{x}z^{-1}\Pi_2^*(A_2-I_n)^{-1},
\end{equation}
 and equations  (\ref{g3'})--(\ref{g5'}) take the form
\begin{equation}       \label{g24}
\big(\Pi_1\big)_t=-(A_1+I_n)^{-1}\Pi_1 z_{t}z^{-1}, \quad \big(\Pi_2^*\big)_x=z_{t}z^{-1}\Pi_2^*(A_2+I_n)^{-1},
\end{equation}
\begin{equation}
  \label{g25}
S_t=(A_1+I_n)^{-1}\Pi_1 z_{t}z^{-1}\Pi_2^*(A_2+I_n)^{-1}.
\end{equation}
Now, let matrices $A_k$, $S(0,0)$, and $\Pi_k(0,0)$
be fixed, assume that (\ref{g18}) holds, and let matrix functions $S$ and $\Pi_k$ satisfy
(\ref{g22})--(\ref{g25}). By Theorem \ref{GBDTrd} we get
\begin{equation}       \label{g26}
\big(w_A\big)_x=\wt Gw_A-w_AG, \quad \big(w_A\big)_t=\wt Fw_A-w_AF.
\end{equation}
Assume additionally that $\det A_k \not= 0$ $(k=1,2)$. Then, by  (\ref{g27'}) and (\ref{g27}) 
the matrix functions $w_A(x,t,0)$ and  $w_A(x,t,0)^{-1}$
are well-defined
in the points of iinvertibility of $S$. 
Taking into account (\ref{g9}) and first equalities in (\ref{g20}) and (\ref{g21}) we have
\begin{equation}       \label{g33}
\wt G(x,t,\la)=-(\la-1)^{-1}\wt q_{11}(x,t), \quad \wt F(x,t,\la)=-(\la+1)^{-1}\wt Q_{11}(x,t).
\end{equation}
It follows from (\ref{g26}) and (\ref{g33}) that
\begin{eqnarray}       \label{g28}&&
\frac{\d}{\d x}w_A(x,t,0)=\wt q_{11}(x,t)w_A(x,t,0)-w_A(x,t,0)q_{11}(x,t), \\
&&  \label{g29}
 \frac{\d}{\d t}w_A(x,t,0)=-\wt Q_{11}(x,t)w_A(x,t,0)+w_A(x,t,0)Q_{11}(x,t).
\end{eqnarray}
Rewrite second relations in  (\ref{g20}) and (\ref{g21}):
\begin{equation}       \label{g30}
z_x=q_{11}z, \quad z_t=-Q_{11}z.
\end{equation}
Put
\begin{equation}       \label{g31}
\wt z(x,t):=w_A(x,t,0)z(x,t).
\end{equation}
From formulas  (\ref{g28})--(\ref{g31}) we derive
$\wt z_x=\wt q_{11}\wt z$ and $\wt z_t=-\wt Q_{11}\wt z$.
As $w_A(x,t,0)$ and $z(x,t)$ are invertible, so $\wt z$ is invertible,
and we get
\begin{equation}       \label{g32}
\wt q_{11}=\wt z_{x}\wt z^{-1}, \quad  \wt Q_{11}=-\wt z_{t}\wt z^{-1}.
\end{equation}
Recall that if formulas (\ref{g20}) and (\ref{g21}) hold, then (\ref{r10}) is equivalent
to (\ref{g19}). The only difference between equalities in (\ref{g20}), (\ref{g21})  and
equalities in (\ref{g33}), (\ref{g32}) is "tilde" in the notations. Hence, in view of (\ref{g33}) and (\ref{g32})
formula (\ref{g18''}) implies that $\wt z$ satisfies  main chiral field equation, that is,
\[
2\wt z_{xt}(x,t)=\wt z_{x}(x,t) \wt z(x,t)^{-1} \wt z_{t}(x,t)+ \wt z_{t}(x,t) \wt z(x,t)^{-1}\wt z_{x}(x,t).
\]
\begin{Cy}\label{CyMCF}   Assume that parameter matrices
satisfy identity (\ref{g18}) and that $\det A_k \not=0$ $(k=1,2)$. 
Let an invertible matrix function $z$ satisfy main chiral field equation
(\ref{g19}) and be two times continuously differentiable.
Then the matrix function 
$\wt z$ given by 
(\ref{g31}) in the points of invertibility of $S$ also satifies main chiral field equation.
\end{Cy}
\end{Ee}
Our next examples deal with the construction
of new (local) solutions of integrable elliptic Sine-Gordon and  sinh-Gordon equations
from the initial solutions.
See, for instance, \cite{BoKi, JaKa} and references therein for some related literature
and auxiliary systems .

\begin{Ee}  \label{EeSG}
Elliptic Sine-Gordon  equation 
\begin{equation}       \label{g36}
v_{tt}+v_{xx}=\sin \, v \quad (v=\ov v)
\end{equation}
is equivalent to the compatibility condition  (\ref{r10})
of the auxiliary systems  (\ref{g15}), where
\begin{eqnarray}       \label{g37}
&& G= \frac{1}{4}\big( i \la \zeta+v_t j-\frac{i}{\la}J\zeta J),
 \quad 
\zeta= \left[
\begin{array}{cc}
0& e^{-iv/2} \\  e^{iv/2} & 0
\end{array}
\right], \\
\label{g38}
&& F=-\frac{1}{4}\big( \la \zeta+v_x j+\frac{1}{\la}J\zeta J),
\end{eqnarray}
and matrices $j$ and $J$ are defined in (\ref{r14}) after putting  $p=1$.
We put also 
\begin{equation}       \label{g39}
A_1=A, \quad A_2=-(A^*)^{-1}, \quad \Pi_1\equiv \Pi,  \quad \Pi_2(0,0)=A^{-1}\Pi(0,0)J.
\end{equation}
Thus, we have three parameter matrices, that is, $n \times n$ matrices $A$ and $S(0,0)$
and an $n \times m$ matrix $\Pi(0,0)$.
We assume that $v$ satisfies (\ref{g36}), that $\det A\not=0$ and $S(0,0)=S(0,0)^*$, and that there is a matrix $U$ such that equalities
\begin{equation}       \label{g40}
\ov A=UA^{-1}U^{-1}, \quad \ov \Pi(0,0)=U\Pi(0,0), \quad \ov S(0,0)=UAS(0,0)A^*U^*
\end{equation}
hold. Here $\ov A$ is the matrix with the entries, which are
complex conjugate to the corresponding entries of $A$.
By  (\ref{g39}) the identity (\ref{g18}), which should be satisfied by the
parameter matrices, takes the form
\begin{equation}       \label{g39'}
AS(0,0)A^*+S(0,0)=\Pi(0,0)J\Pi(0,0)^*.
\end{equation}
Compare (\ref{g16}) and (\ref{g17}) with (\ref{g37}) and (\ref{g38}), respectively,
to see that
\[
r=1, \, \, q_1=(-i/4)\zeta, \,\, q_0=(-v_t/4) j; \,\, l=r_1=1, \,\, c_1=0, \,\, q_{11}=(i/4)J\zeta J;
\]
\[
R=1, \, \, Q_1=(1/4)\zeta, \,\, Q_0=(v_x/4) j; \,\, L=R_1=1, \,\, C_1=0, \,\, Q_{11}=(1/4)J\zeta J.
\]
Thus, in view of   (\ref{g39}) equations  (\ref{g3}) and  (\ref{g3'})
take the form
\begin{equation}       \label{g41}
\Pi_x=\frac{1}{4}\big(-iA\Pi \zeta -v_t\Pi j+iA^{-1}\Pi J\zeta J\big), \quad
\Pi_t=\frac{1}{4}\big(A\Pi \zeta +v_x\Pi j+A^{-1}\Pi J\zeta J\big).
\end{equation}
As $\ov \zeta = J\zeta J$ and $\ov A=UA^{-1}U^{-1}$, one can see that both
$\Pi$ and $z=\ov U \ov \Pi$ satisfy (\ref{g41}).  According to (\ref{g40}) $\Pi(0,0)=\ov U \ov \Pi(0,0))$, 
and so we derive
\begin{equation}       \label{g41'}
 \Pi(x,t) \equiv  \ov U \ov \Pi(x,t) .
\end{equation}
Equations (\ref{g4}) and  (\ref{g4'}), which define $\Pi_2^*$, take the form
\begin{eqnarray}       \label{g42}
&&(\Pi_2^*)_x=\frac{1}{4}\big(i\zeta\Pi_2^*A_2 +v_t j\Pi_2^*-i J\zeta J\Pi_2^*A_2^{-1}\big), \\
 \label{g43} &&
(\Pi_2^*)_t=\frac{1}{4}\big(-\zeta\Pi_2^*A_2 -v_x j\Pi_2^*-J\zeta J\Pi_2^*A_2^{-1}\big).
\end{eqnarray}
As $\zeta=\zeta^*$ and $A_2=-(A^*)^{-1}$, it follows from (\ref{g41})
that 
\[
z(x,t)=A^{-1}\Pi(x,t)J
\]
satisfies  equations (\ref{g42}) and  (\ref{g43}) for
$\Pi_2$.  Moreover, we have $\Pi_2(0,0)= A^{-1}\Pi(0,0)J=z(0,0)$.
In other words we have
\begin{equation}       \label{g44}
\Pi_2(x,t) \equiv A^{-1}\Pi(x,t)J,
\end{equation}
and identity  (\ref{g18'}) takes the form
\[
AS(x,t)A^*+S(x,t)=\Pi(x,t)J\Pi(x,t)^*.
\]
By   (\ref{g5}),  (\ref{g5'}),   (\ref{g44}), and by the second equality in   (\ref{g39}) the relations
\begin{equation}       \label{g45}
S_x=\frac{i}{4}\big(A^{-1}\Pi J\zeta \Pi^*-\Pi \zeta J\Pi^*(A^*)^{-1}\big),
\,\,  S_t=\frac{1}{4}\big(A^{-1}\Pi J\zeta \Pi^*+\Pi \zeta J\Pi^*(A^*)^{-1}\big)
\end{equation}
hold. Formulas (\ref{g40}),  (\ref{g41'}), and (\ref{g45}) imply
$\ov S_x=UAS_xA^*U^*$, $\ov S_t=UAS_tA^*U^*$ and finally
\begin{equation}       \label{g46}
\ov S \equiv UASA^*U^*.
\end{equation}
It follows from (\ref{g12}), (\ref{g41'}), (\ref{g44}), and (\ref{g46}) that
\begin{equation}       \label{g47}
X_{-1}=J\Pi^*(A^*)^{-1}S^{-1}A^{-1}\Pi=J\Pi^*U^*(UASA^*U^*)^{-1}U\Pi=J\ov{\big(\Pi^* S^{-1}\Pi\big)}.
\end{equation}
According to  (\ref{g12}),  (\ref{g27'}), and (\ref{g27}) we have 
\begin{equation}       \label{g48}
Z(x,t):=w_A(x,t,0)=I_2-X_{-1}, \quad Z(x,t)^{-1}=w_A(x,t,0)^{-1}=I_2+Y_{-1}.
\end{equation}
Moreover, in view of (\ref{g12}), (\ref{g44}), and equality $\Pi^* S^{-1}\Pi=(\Pi^* S^{-1}\Pi)^*$ we get
\begin{equation}       \label{g49}
Y_{-1}=-J\Pi^* S^{-1}\Pi, \quad -\Pi^* S^{-1}\Pi=\left[
\begin{array}{cc}
a & b \\ \ov b & d
\end{array}
\right], \quad a=\ov a, \quad d=\ov d.
\end{equation}
Using  (\ref{g47})--(\ref{g49}) we derive
\begin{equation}       \label{g50}
Z\, Z^{-1}=\left[
\begin{array}{cc}
1+b & d \\  a & 1+\ov b
\end{array}
\right]
\left[
\begin{array}{cc}
1+\ov b & d \\  a & 1+ b
\end{array}
\right]=I_2.
\end{equation}
If $1+b \not=0$, formula (\ref{g50})  implies
\begin{equation}       \label{g51}
a=d=0, \quad  |1+b|=1, \quad Z=\diag\{1+b,\, 1+\ov b\}.
\end{equation}
Put
\begin{equation}       \label{g52}
\wh v=v+2 \arg(1+b), \quad \wh u=Z^{-\frac{1}{2}}\wt u, \quad \wh G= \wh u_x  \wh u^{-1},
\quad \wh F= \wh u_t  \wh u^{-1},
\end{equation}
where $\wt u(x,t,\la)=w_A(x,t,\la)u(x,t,\la)$.  In a way, which is similar to the proof of
(\ref{g18''}) in Theorem \ref{ZQE}, we derive from (\ref{g52}) that 
\begin{equation}       \label{g53}
\wh G_t- \wh F_x+[\wh G, \wh F]=0.
\end{equation}
Moreover, one can see that $\wh G$ and $\wh F$ have the form (\ref{g37})
and (\ref{g38}), respectively, after one substitutes $\wh v$ instead of $v$
into the right-hand sides of (\ref{g37})
and (\ref{g38}). Therefore, formula (\ref{g53}) implies that $\wh v$ satisfies
(\ref{g36}).
\begin{Cy}\label{eSG} Let an integer  $n>0$ and matrices
$A$ ($\det A\not=0$), $\Pi(0,0)$, and $S(0,0)=S(0,0)^*$ be fixed and satisfy conditions 
(\ref{g40}) and (\ref{g39'}). 
Let $v$ satisfy elliptic Sine-Gordon equation
and be two times continuously differentiable. 
Then in the points, where
$\det S \not=0$ and $1+b\not=0$, the function $\wh v$ given by (\ref{g52})
satisfies elliptic Sine-Gordon equation too.

If  $\s(A) \cap \s\big((-A^*)^{-1}\big)=\emptyset$, then the last equality in 
 (\ref{g40}) follows from  (\ref{g39'}) and from the first two equalities in  (\ref{g40}).
\end{Cy}
\end{Ee}
\begin{Rk}\label{Rk-1}
Using considerations from (\ref{g48}) one easily shows  that a general equality  $(I_m-X_{-1})(I_m+Y_{-1})=I_m$
is true.
\end{Rk}
Elliptic sinh-Gordon  equation 
\begin{equation}       \label{g54}
v_{tt}+v_{xx}=\sinh \, v \quad (v=\ov v)
\end{equation}
is equivalent to the compatibility condition  (\ref{r10})
of the auxiliary systems  (\ref{g15}), where
\begin{eqnarray}       \label{g55}
&& G=- \frac{1}{4}\big(  \la \zeta-iv_t j+\frac{1}{\la}\zeta^*),
 \quad 
\zeta= \left[
\begin{array}{cc}
0& e^{-v/2} \\  e^{v/2} & 0
\end{array}
\right], \\
\label{g56}
&& F=-\frac{1}{4}\big( -i\la \zeta+iv_x j+\frac{i}{\la}\zeta^*).
\end{eqnarray} 
Put
\begin{equation}       \label{g57}
A_1=A, \quad A_2=-(A^*)^{-1}, \quad \Pi_1\equiv \Pi,  \quad \Pi_2(0,0)=A^{-1}\Pi(0,0).
\end{equation}
Here we assume that
\begin{equation}       \label{g58}
\det\, A \not=0, \quad S(0,0)=S(0,0)^*,
\end{equation}
and that there is a matrix $U$ such that
\begin{equation}       \label{g59}
\ov A=UA^{-1}U^{-1}, \quad \ov \Pi(0,0)=U\Pi(0,0)J, \quad \ov S(0,0)=UAS(0,0)A^*U^*.
\end{equation}
Now, the identity (\ref{g18}) takes the form
\begin{equation}       \label{g60}
AS(0,0)A^*+S(0,0)=\Pi(0,0)\Pi(0,0)^*.
\end{equation}
Taking into account (\ref{g55}) and (\ref{g56}) introduce $\Pi$ by the equations
\begin{equation}       \label{g61}
\Pi_x=\frac{1}{4}\big(A\Pi \zeta -iv_t\Pi j+A^{-1}\Pi \zeta^* \big), \quad
\Pi_t=\frac{1}{4}\big(-iA\Pi \zeta +iv_x\Pi j+iA^{-1}\Pi \zeta^* \big).
\end{equation}
It is easy to see that $\Pi_2 \equiv A^{-1}\Pi$, and so formulas  (\ref{g5}) and (\ref{g5'})
take the form
\begin{equation}       \label{g62}
S_x=\frac{1}{4}\big(\Pi \zeta \Pi^*(A^*)^{-1}+A^{-1}\Pi \zeta^*\Pi^*\big),
\,\,  S_t=\frac{i}{4}\big(A^{-1}\Pi \zeta^*\Pi^*- \Pi \zeta \Pi^*(A^*)^{-1}\big).
\end{equation}
The matrix $Z=I_2-X_{-1}$ is again a diagonal matrix, and we have
\begin{equation}       \label{g63}
 Z=I_2-\Pi^*(A^*)^{-1}S^{-1}A^{-1}\Pi=\diag\{Z_{11},\, Z_{11}^{-1}\}, \quad Z_{11}=\ov{Z_{11}} .
\end{equation}
(compare with formula (\ref{g51}), the proof of (\ref{g63}) is similar).
The following corollary is proved in a quite similar way to the Corollary \ref{eSG}.
\begin{Cy}\label{eshG} Let an integer  $n>0$ and matrices
$A$, $\Pi(0,0)$, and $S(0,0)$ be fixed and satisfy conditions 
(\ref{g58})--(\ref{g60}). 
Let $v$ satisfy elliptic sinh-Gordon equation
and be two times continuously differentiable. 
Then in the points, where
$\det S \not=0$ and $Z_{11}\not=0$, 
the function 
\begin{equation}       \label{g64}
\wh v=v+2\ln |Z_{11}|
\end{equation}
satisfies elliptic Sine-Gordon equation too. Here $Z_{11}$ is given by 
(\ref{g61})--(\ref{g63}).
\end{Cy}

\section{GBDT for  radial Dirac equation} \label{DirRad}
\setcounter{equation}{0}
Radial Dirac equation has the form
\begin{equation}       \label{1.0}
\Big(-i\s_2\frac{d}{dx}+\frac{\k}{x}\s_1+V(x)\Big)u=\la u
\quad
(x > 0),
\end{equation}
or equivalently
\begin{equation}       \label{1.1}
\Big(\frac{d}{dx}+\la q_1 + q_0(x)\Big)u(x, \lambda )=0 \quad
(x > 0),
\end{equation}
where $\kappa$ is integer, $\s_i$  are Pauli matrices,
\begin{equation}   \label{1.2}
\s_1 = \left[
\begin{array}{cc}
0 & 1 \\ 1 & 0
\end{array}
\right], \quad \s_2 = \left[
\begin{array}{cc}
0 & -i \\ i & 0
\end{array}
\right], \quad \s_3 = \left[
\begin{array}{cc}
1 & 0 \\ 0 & -1
\end{array}
\right],
 \end{equation}
\begin{equation}   \label{1.2'}
V(x)=v_e(x)I_2+v_a(x)\s_1+v_s(x)\s_3, 
 \end{equation} 
$v_e$, $v_a$, and $v_s$ are real-valued functions, which are  locally integrable on $[0, l)$  (here and further $l \in \BR_+$),
\begin{equation}       \label{1.3}
q_1=-\breve J, \quad \breve J:=(-i\s_2)^{-1}= \left[
\begin{array}{cc}
0 & 1 \\ -1 & 0
\end{array}
\right], 
\end{equation}
\begin{equation}
 \label{1.3'}  q_0(x)=v_*(x) \s_3+ \breve J \big(v_e(x)I_2+v_s(x)\s_3\big), \quad v_*(x)=\frac{\k}{x}+v_a(x).
\end{equation}
Here $v_e$, $v_a$, and $v_s$ represent the electrostatic potential, the anomalous magnetic moment,
and the sum of the mass and the scalar potential, respectively.

If we put $p=1$ in (\ref{r14}) we have $\s_1=J$ and $\s_3=j$.
Recall that  Dirac-type system of the form (\ref{0x1.-4}), (\ref{0x1.-2})
was treated in Subsection \ref{subGauge}. The radial Dirac equation
differs from the Dirac-type system. Its structure (as well as the structure of $V$
in this section) is somewhat different and it usually has singularity at $x=0$,
which is of interest from the physical point of view. We consider equation  (\ref{1.0}) 
 independently from the results in  the Subsection \ref{subGauge}.

The double commutation method was applied to (\ref{1.0}) in an
interesting paper by G. Teschl \cite{T}. By this method  S. Albeverio, R. Hryniv, and Ya. Mykytyuk
\cite{AHM} proved that $1$ is added to $\vk$, when an eigenvalue
is removed and that $1$ is subtracted from $\vk$, when an eigenvalue is inserted.
(In this section we actively use some of the results from \cite{AHM}.)
We apply GBDT to the radial Dirac equation (\ref{1.0}). In particular,
we construct  explicitly potentials and fundamental solutions
for the  equation (\ref{1.0}) with $\vk >0$ starting from the
trivial equation (i.e., from $q_0=0$).  The case formally corresponds
to the removal of the eigenvalues, the iterated double commutation formulas
for the insertion of the eigenvalues are given in \cite{T}.

Fundamental solution $u$ in this section is a non-degenerate $2 \times 2$ solution of
(\ref{1.0}),  we do not require $u$ to be  normalized at $x=0$.

\subsection{Main result} \label{MR}
The following procedure to construct explicit solutions of  
the radial Dirac equation is an immediate corollary of   Theorem \ref{Tm1}
from Subsection \ref{Au}.
\begin{Tm}  \label{Expl}
To construct a class of  explicit solutions of equation (\ref{1.0}) with some fixed integer
$\k$, fix  an integer $m>0$ and  $m \times m$  matrices ${\cal A}_1$
and ${\cal S}_1>0$. Fix also a $\vk \times \vk$ lower triangular  matrix ${\cal A}_2$ and an $(m+\vk)\times 2$
matrix $\Pi(0)={\mathrm{col}}[\Psi_1(0) \quad \Psi_2(0)]$, where $\vk=|\k|$,
col means column and $\Psi_1$ $(\Psi_2)$ is an $m \times 2$ upper $(\vk \times 2$ lower$)$ block of $\Pi$.
It is required that
\begin{equation}       \label{e1}
{\cal A}_1{\cal S}_1-{\cal S}_1{\cal A}_1^*=\Psi_1(0)\breve J\Psi_1(0)^*, \quad  \Psi_2(0)\breve J\Psi_2(0)^*=0.
\end{equation}
Moreover, for  $h: =[1 \quad 0 \quad 0 \, \ldots \, 0]\Psi_2$  we assume that
$ h(0)=c[0\quad 1]$, when $\k$ is positive and odd or negative and even and 
that $h(0)=c[1\quad 0]$, when $\k$ is positive and even or negative and odd $(c\not=0)$.

Introduce matrix $A$
\begin{equation}       \label{e2}
A:=\left[
\begin{array}{lr}
{\cal A}_1& 0 \\ R & {\cal A}_2
\end{array}
\right], \quad R=\Psi_2(0)\breve J\Psi_1(0)^*{\cal S}_1   ^{-1}, 
\end{equation}
and vectors $ \t_1, \, \t_2 \in \BC^{m+\vk}$
\begin{equation}       \label{e3}
[\t_1 \quad \t_2]:=\Pi(0)\breve K, \quad \breve K:=\frac{1}{\sqrt{2}}\left[\begin{array}{lr}
1 & 1\\ -i & i
\end{array}
\right].
\end{equation}
Now, put
\begin{eqnarray}       \label{e3'}
&& \Pi(x):=[e^{ixA}\t_1 \quad e^{-ixA}\t_2]\breve K^*, \\
 \label{e3d}
&& S(x):=\left[
\begin{array}{lr}
{\cal S}_1& 0 \\ 0 & 0
\end{array}
\right]+\int_0^x \Big(e^{itA}\t_1 \t_1^*e^{-itA^*}+e^{-itA}\t_2 \t_2^*e^{itA^*}\Big)dt.
\end{eqnarray}
Then, if $S(x)>0$ for $x>0$,  the potential
\begin{equation}       \label{e4}
\wt q_0(x):=\breve JX(x)\breve J^*-X(x), \quad X=\{X_{ij}\}_{i,j=1}^2:=\Pi^*S^{-1}\Pi
\end{equation}
admits representation 
\begin{equation}       \label{e4'}
\wt q_0(x)=\frac{\k}{x}\s_3+\Up(x)=\frac{\k}{x}\s_3+\wt v_a(x)\s_3-\wt v_s(x)\s_1,
\end{equation}
where  $\Up$, $\wt v_a$, and  $\wt v_s$ are bounded in the neighborhood of zero, and
\begin{equation}       \label{e4''}
 \wt v_a=X_{22}-X_{11}-\frac{\k}{x},
\quad  \wt v_s=X_{12}+X_{21}.
\end{equation}
The fundamental solution $\wt u(x, \la) $ of  system 
\begin{equation}       \label{1.1'}
\Big(\frac{d}{dx}+\la q_1 + \wt q_0(x)\Big)\wt u(x, \lambda )=0 \quad
(x > 0),
\end{equation}
is given by the formula
\begin{equation}       \label{e5}
\wt u(x,\la) =w_A(x,\la) \breve K\exp({-i\la x \s_3}),
\end{equation}
where 
\begin{equation}       \label{2.1}
w_A(\lambda )=I_{2}-\breve J \Pi^*S^{-1}\big(A- \lambda I_n
\big)^{-1} \Pi.
\end{equation}
\end{Tm}
By the neighbourhood of zero we  mean the neighbourhood of the form $(0, \, \ve)$
or  $[0, \, \ve)$. Note that if $S(x)>0$ for all $x \in (0, \, \ve)$, then $S(x)>0$ for  all $x \in \BR_+$,
the formulas in Theorem \ref{Expl} are well-defined on $\BR_+$, and $\wt q_0$
is infinitely differentiable on $\BR_+$.
\subsection{Superposition of Darboux transformations} \label{Au}

In this subsection we apply GBDT to the radial Dirac equation.
We factorize also the  Darboux matrix, so that GBDT can be cosidered
as a superposition of two other GBDTs. Finally, we formulate a general 
Theorem \ref{Tm1}.

In view of  (\ref{1.1}),  (\ref{1.3}), and (\ref{1.3'})  we get
\begin{equation}       \label{r16}
r=1, \quad q_1=-\breve J; \quad q_k^*=-\breve Jq_k\breve J^{-1}, \quad k=0, \, 1, 
\quad \breve J^{-1}=\breve J^*=-\breve J.
\end{equation}
Here we  fix $n>0$ and parameter matrices  $A_1=A$, $\Pi(x_0)$, and $S(x_0)=S(x_0)^*$.
Formula (\ref{g3}) for $\Pi_1=\Pi$ takes the form
\begin{equation}       \label{2.2}
 \Pi_x=A\Pi q_1+\Pi q_0. 
\end{equation}
Putting
\begin{equation}       \label{r17}
A_2=A^*, \quad \Pi_2^*=\breve J \Pi^*
\end{equation}
and taking into account  equation (\ref{2.2}) and the third equalities in  (\ref{r16}),
we see that equation  (\ref{g4}) for $\Pi_2^*$ is satisfied.
By (\ref{2.2}) and (\ref{r17}) formula (\ref{g5})  takes the form
\begin{equation}       \label{2.2'}
 S_x=\Pi\Pi^*.
\end{equation}
Thus, we have $S(x)=S(x)^*$. It is required that
\begin{equation}       \label{2.4}
AS(x_0)-S(x_0)A^*=\Pi(x_0)\breve J\Pi(x_0)^*,
\end{equation}
$0<x_0<l. $ The identity
\begin{equation}       \label{2.3}
AS(x)-S(x)A^*=\Pi(x)\breve J\Pi(x)^*
\end{equation}
follows from (\ref{2.2}), (\ref{2.2'}), and (\ref{2.4}).
\begin{Rk}\label{Rk1}
If  $\Pi$ can be continuously extended to 
$\Pi(0)$, then $S$ can also be continuously extended to $S(0)$. So, in that case we consider
$\Pi$ and $S$  defined on $[0, \, l)$, and for  (\ref{2.2}) to be true it suffices
that (\ref{2.4}) holds for $x_0=0$. 
\end{Rk}
The next  Corollary of Theorem \ref{GBDTrd} and formula (\ref{g14})  is immediate
\begin{Cy}  \label{Pn1}
Let relations  (\ref{2.2}), (\ref{2.2'}), and  (\ref{2.4}) hold. Then, in the points of the invertibility of $S(x)$ $(x>0)$, we have
\begin{equation}       \label{2.6}
\Big(S(x)^{-1}\Pi(x)\Big)_x(x)=A^*S(x)^{-1}\Pi(x) q_1(x)+S(x)^{-1}\Pi(x)\wt q_0(x),
\end{equation}
where
\begin{equation}       \label{2.7}
\wt q_0:=q_0+\breve JX\breve J^*-X, \quad X:=\breve J^{-1}X_0=\Pi^*S^{-1}\Pi.
\end{equation}
Moreover, the matrix function $w_A$ given by (\ref{2.1})  satisfies the equation
\begin{equation}       \label{2.8}
\frac{d}{dx}w_A(x, \lambda )=\wt G(x,\la) w_A(x, \lambda )-w_A(x, \lambda )G(x,\la) , 
\end{equation}
\begin{equation}       \label{2.9}
\wt G(x,\la) =-\la q_1-\wt q_0(r), \quad G(x,\la) =-\la q_1- q_0(x).
\end{equation}
In other words, we have
\begin{equation}       \label{2.10}
\frac{d}{dx}\wt u(x, \la) =\wt G(x,\la) \wt u(x, \la) , 
\end{equation}
where 
\begin{equation}       \label{2.11}
\wt u(x, \la) =w_A(x,\la)  u(x, \la) .
\end{equation}
\end{Cy}
Recall that by Remark \ref{Rkr} $\wt q_1=q_1$. This explains why we have
the coefficient $-q_1$ in the expression for $\wt G$ in (\ref{2.9}).
Equation (\ref{2.10})
follows from (\ref{1.1}) and  (\ref{2.8}). This equation is a GBDT transformation
of the equation  (\ref{1.1}) and has the same structure. Namely, to get 
$\wt q_0$ instead of $q_0$
one substitutes into (\ref{1.3'})  $\wt v_*$ and $\wt v_s$  instead of  $v_*$ and $v_s$, respectively.
That is, by  (\ref{1.3'}) and (\ref{2.7}) we have
\begin{equation}
 \label{r18}  \wt q_0(x)=\wt v_e(x)\breve J+v_*(x) \s_3- \wt v_s(x)\s_1,
\end{equation}
where $\wt v_e(x)=v_e(x),$
\begin{equation}       \label{2.12}
  \wt v_*(x)=v_*(x)+X_{22}(x)-X_{11}(x), 
\quad  \wt v_s(x)=v_s(x)+X_{12}(x)+X_{21}(x),
\end{equation}
and $X_{kj}$ are the entries of $X$. 

\begin{Rk}\label{Rk3}
When $n=1$ and $A \in \BR$ our transformation coincides
with the double commutation transformation for Dirac equations
treated in \cite{AHM, T}. When $A \not\in \BR$ our transformation
somewhat differs from the transformation in \cite{AHM}, because
the transformation in \cite{AHM} uses the equivalent (for $n=1$) of $\Pi$
and of the transposition of $\Pi$, whereas we use here $\Pi$ and $\Pi^*$,
so that $\wt v_a$ and $\wt v_s$ are real valued.
\end{Rk}
We will be interested in the transformation of the equation (\ref{1.1})
with $\kappa =0$ into equation with integer nonzero $\kappa$.
\begin{Ee} \label{Ee3} Let $n=1$, $q_0 \equiv 0$, $S(0)=0$.  Notice that we have
\begin{equation}       \label{2.13}
q_1=-\breve J=i\breve K\s_3 \breve K^*,  \quad \breve K:=\frac{1}{\sqrt{2}}\left[\begin{array}{lr}
1 & 1\\ -i & i
\end{array}
\right], \quad \breve K^*=\breve K^{-1},
\end{equation}
and put $(\Pi \breve K)(0)=[1 \quad \a]$. Taking into account $q_0=0$ and  (\ref{2.13}), we  rewrite (\ref{2.2}) 
and (\ref{2.2'})
in the form
\begin{equation}       \label{2.14}
(\Pi \breve K)_x=iA\Pi \breve K \s_1, \quad S_x=(\Pi \breve K)(\Pi \breve K)^*
\end{equation}
 It is immediate from  (\ref{2.14}) that
\begin{equation}       \label{2.15}
\Pi(x)\breve K=[e^{ixA} \quad \a e^{-ixA}].
\end{equation}
Recall that $S(0)=0$. Thus, when $A \not= \ov A$, we require additionally $|\a|=1$  so that  (\ref{2.3}) holds at $x=0$. Hence, formula
(\ref{2.3}) is true  for all $x>0$ (see Remark \ref{Rk1}). Using (\ref{2.3}), (\ref{2.14}) and (\ref{2.15}), we get
\begin{eqnarray}       \nonumber
&& S(x)=i(A-\ov A)^{-1}\Big( e^{ix(\ov A-A)}  -e^{ix(A-\ov A)} \Big) \quad {\mathrm{for}} \quad A \not= \ov A, \\
&& S(x)=(1+|\a|^2)x \quad {\mathrm{for}} \quad A\in \BR. \label{2.16}
\end{eqnarray}
By  (\ref{2.12}), (\ref{2.15}) and (\ref{2.16}) we have for  $A \not= \ov A$ the equality
\begin{equation}       \label{2.17}
\wt v_*(x)=i( A-\ov  A)\Big(\a e^{-ix(A+\ov A)}+\ov \a e^{ix(A+\ov A)}\Big)\Big(e^{ix(\ov A-A)}  -e^{ix(A-\ov A)}\Big)^{-1}.
\end{equation}
In a similar way we have
\begin{equation}       \label{2.18}
\wt v_*(x)=-\Big(\a e^{-2ixA}+\ov \a e^{2ixA}\Big)\Big((1+|\a|^2)x\Big)^{-1} \quad {\mathrm{for}} \quad A \in \BR.
\end{equation}
From (\ref{2.17}) and (\ref{2.18}) we derive that
\begin{equation}       \label{2.19}
\wt v_*(x)=\frac{\wt \k}{x}+\wt v_a, \quad \wt \kappa =-\frac{\a+\ov \a}{1+|\a|^2},
\end{equation}
where $\wt v_a$ is continuous on $[0, \, \infty)$. When $\a=\pm 1$, we obtain $\wt \kappa =\mp 1$, i.e., $\kappa$ is integer.
\end{Ee}
To study GBDT we will need, in particular, to split it into a superposition of several transformations,
and our next result is dedicated to this procedure. Let $A$ be a block lower triangular matrix:
\begin{equation}       \label{2.20}
A= \left[
\begin{array}{lr}
A_{11}& 0\\ A_{21} & A_{22}
\end{array}
\right], \quad S= \left[
\begin{array}{lr}
S_{11}& S_{12}\\ S_{21} & S_{22}
\end{array}
\right],
\end{equation}
where $A_{11}$ and $S_{11}$ are $n_1 \times n_1$ matrices,
$A_{22}$ and $S_{22}$ are $n_2 \times n_2$ matrices, $n=n_1+n_2$.
In our further considerations we fix some value of $x$ and omit temporarily
the variable $x$ in the notations.
Assume  $\det S\not=0$ and $\det S_{11}\not=0$,  and denote by
$T_{22}$  the $n_2 \times n_2$ right  lower block of $T=S^{-1}$.
The invertibility of
$T_{22}$ follows from the invertibility of $S_{11}$ and $S$. One can  check directly that
\begin{equation}  \label{2.21}
T=
\left[
     \begin{array}{cc}
             S_{11}^{-1}+    S_{11}^{-1}S_{12} \, T_{22} \, S_{21}  S_{11}^{-1}
          & -  S_{11}^{-1} S_{12}  T_{22}
     \\
            - T_{22}  S_{21}  S_{11}^{-1}
          & T_{22}
     \end{array}
\right], \quad T_{22}^{-1}= S_{22}-S_{21}S_{11}^{-1}S_{12}.
\end{equation}
As $A$ is a block lower triangular matrix and  $\det S\not=0$, 
$\det S_{11}\not=0$,  so $w_A$ admits factorisation \cite{SaL1, SaL3}
\begin{equation}       \label{2.22}
w_A( \la) =w_2( \la) w_1( \la) ,
\end{equation}
\begin{eqnarray} \label{2.23}
&& w_1(\la) =I_2-\breve J\Pi^*P_1^*S_{11}^{-1}\big(A_{11}-\la I_{n_1})^{-1}P_1 \Pi , \\ \label{2.23'}
&& \nonumber w_2(\la) =I-\breve J\Pi^*S^{-1}P_2^*\big(A_{22}-\la I_{n_2})^{-1}T_{22}^{-1}P_2S^{-1} \Pi,
\end{eqnarray}
where 
\begin{equation}       \label{2.24}
P_1=\left[
     \begin{array}{cc}
I_{n_1} & 0
     \end{array}
\right], \quad P_2=\left[
     \begin{array}{cc} 0 &
I_{n_2}
     \end{array}
\right].
\end{equation}
One can easily see that by  (\ref{2.20}) we have $P_1 A=A_{11}P_1$ and $AP_2^*=P_2^*A_{22}$.
Hence,
the identity
\begin{equation}       \label{2.25}
A_{11}S_{11}-S_{11}A_{11}^*= \pi_1 \breve J\pi_1^*, \quad \pi_1:=P_1 \Pi
\end{equation}
follows from  (\ref{2.3}). Now, compare (\ref{2.1}) and the first relation in (\ref{2.23}) to see that  $w_1=w_{A_{11}}$.
That is, we get $w_1$ after substitution  of $A_{11}$ instead of $A$,
of $S_{11}$ 
   instead of $S$ and of $\pi_1$  instead of  $\Pi$ into (\ref{2.1}). 
Identity  (\ref{2.25}) is the equivalent of  (\ref{2.3}) written for the new transfer matrix function
$w_{A_{11}}$.

Moreover,
from (\ref{2.3})  it follows that $TA-A^*T=T\Pi \breve J \Pi^* T$. Therefore, we derive
$T_{22}A_{22}-A_{22}^*T_{22}=P_2T\Pi \breve J \Pi^* TP_2^*$ or, equivalently
\begin{equation}       \label{2.26}
A_{22}T_{22}^{-1}-T_{22}^{-1}A^*_{22}=\pi_2 \breve J \pi_2^*, \quad \pi_2=T_{22}^{-1}P_2S^{-1}\Pi.
\end{equation}
Compare  (\ref{2.1})  and (\ref{2.23})  to   see that $w_2=w_{A_{22}}$, where $T_{22}^{-1}$ is substituted into
 (\ref{2.1})  instead of $S$ and $\pi_2$ is substituted instead of  $\Pi$. 

\begin{Rk} \label{sup}
Corollary \ref{Pn1} defines by formulae  (\ref{2.7}) and (\ref{2.11})  GBDT
transformation of the coefficient $q_0$ and  solutions $u$ of the radial Dirac
equation (\ref{1.1}).  The next proposition shows  that the GBDT generated by 
the parameter matrices $A$, $S(x_0)$ and $\Pi(x_0)$  can be treated as a superposition
of two transformations, i.e., GBDTs  generated by $A_{11}$, $S_{11}(x_0)$, $\pi_1(x_0)$ 
and $A_{22}$, $T_{22}(x_0)^{-1}$, $\pi_1(x_0)$, respectively.
In particular,  our next 
proposition shows that $\pi_1$ and $S_{11}$  satisfy analogs of  (\ref{2.2}) and (\ref{2.2'}), respectively.
 The same is true
for  $T_{22}^{-1}$ and $\pi_2$. According to Corollary \ref{Pn1} it means that
$w_1(x,\la) $ and $w_2(x,\la) $ are Darboux matrices. 
\end{Rk}  
\begin{Pn}  \label{fact} Let relations  (\ref{2.2}), (\ref{2.2'}), and  (\ref{2.4}) be valid and let
$A$ be a block lower triangular matrix. Then, in the points of the invertibility of $S(x)$ and $S_{11}(x)$ $(x>0)$, we have
\begin{equation}       \label{2.27}
\big(S_{11}\big)_x=\pi_1\pi_1^*, \quad (\pi_1)_x=A_{11}\pi_1 q_1+\pi_1 q_0.
\end{equation}
\begin{equation}       \label{2.28}
\Big(T_{22}^{-1}\Big)_x=\pi_2\pi_2^*, \quad (\pi_2)_x=A_{22}\pi_2 q_1+\pi_2 \wh q_0,
\end{equation}
where $\pi_1$ and $\pi_2$ are given by  (\ref{2.25}) and  (\ref{2.26}), respectively, and
\begin{equation}       \label{2.29}
\wh q_0=q_0+\breve J\wh X \breve J^*- \wh X, \quad \wh X=\pi_1^*S_{11}^{-1}\pi_1.
\end{equation}
\end{Pn}
\begin{proof}.
Multiply from the left both sides of   (\ref{2.2}) by $P_1$ and use $P_1 A=A_{11}P_1$ 
to get  the second relation in (\ref{2.27}).
Multiply by $P_1$ from the left and by $P_1^*$ from the right both sides of
 (\ref{2.2'})  to get the  first relation in (\ref{2.27}).
In view of 
 (\ref{2.2'}) and of the definition of $\pi_2$ in (\ref{2.26}), we have also the  first relation in (\ref{2.28}):
\begin{eqnarray}        \nonumber
&&\frac{d}{dx}T_{22}^{-1}=-T_{22}^{-1}\Big(\frac{d}{dx}T_{22}\Big)T_{22}^{-1}
=-T_{22}^{-1}P_2\Big(\frac{d}{dx}T\Big)P_2^*T_{22}^{-1}
\\ 
&&=T_{22}^{-1}P_2S^{-1}\Pi \Pi^*S^{-1}P_2^*T_{22}^{-1}=\pi_2\pi_2^*.  \label{2.30}
\end{eqnarray}
Now, use (\ref{2.6}),  (\ref{2.30}), and equality $P_2A^*=A_{22}^*P_2$ to differentiate $\pi_2$:
\begin{eqnarray}     \nonumber  
&&(\pi_2)_x=\pi_2\pi_2^*P_2S^{-1}\Pi+T_{22}^{-1}P_2\big(A^*S^{-1}\Pi q_1+S^{-1}\Pi\wt q_0 \big) \\
\label{2.31}
&&=\pi_2\pi_2^*P_2S^{-1}\Pi+T_{22}^{-1}P_2S^{-1}\Pi\wt q_0+T_{22}^{-1}A_{22}^*P_2S^{-1}\Pi q_1.
\end{eqnarray}
Taking into account both relations in (\ref{2.26}), both relations in (\ref{2.7}), and $q_1=\breve J^*$,  rewrite (\ref{2.31}) in the form
\begin{eqnarray}      \nonumber 
&&(\pi_2)_x=\pi_2\pi_2^*P_2S^{-1}\Pi+\pi_2 \wt q_0+A_{22}T_{22}^{-1}P_2S^{-1}\Pi q_1
\\ &&
-\pi_2\breve J\pi_2^*P_2S^{-1}\Pi q_1=A_{22}\pi_2 q_1+\pi_2\big(q_0+\breve J\Pi^*S^{-1}\Pi \breve J^*\nonumber 
\\ &&
-\Pi^*S^{-1}\Pi +\pi_2^*P_2S^{-1}\Pi -\breve J\pi_2^*P_2S^{-1}\Pi \breve J^*
\big).  \label{2.32}
\end{eqnarray}
Finally, note that according to the definition of $\pi_2^*$ and formula (\ref{2.21})  we get
\begin{equation}       \label{2.33}
\Pi^*-\pi_2^*P_2=\Pi^*\left(I_n+\left[
     \begin{array}{c}
S_{11}^{-1}S_{12} \\ -I_{n_2}
     \end{array}
\right]P_2\right)=\pi_1^*S_{11}^{-1}[S_{11} \quad S_{12}].
\end{equation}
It follows that
\begin{equation}       \label{2.34}
\big(\Pi^*-\pi_2^*P_2\big)S^{-1}\Pi=\pi_1^*S_{11}^{-1}\pi_1
\end{equation}
Substitute (\ref{2.34}) into (\ref{2.32}) to derive
\begin{equation}       \label{2.35}
(\pi_2)_x=A_{22}\pi_2 q_1+\pi_2\big(q_0+\breve J\pi_1^*S_{11}^{-1}\pi_1\breve J^*-\pi_1^*S_{11}^{-1}\pi_1\big).
\end{equation}
Formulas (\ref{2.29}) and (\ref{2.35}) imply the second relation in (\ref{2.28}).
\end{proof}
Proposition  \ref{fact} deals with the $S$-nodes, which appear in the process of factorization
(\ref{2.22}) of  the transfer matrix function $w_A$. An inverse in a certain sense
result  is given in the next proposition.
\begin{Pn}  \label{syn}
Let $w_{{\cal A}_j}$  $(j=1,2)$ be  transfer matrix functions of the form (\ref{2.1})
corresponding to the $S$-nodes
${\cal A}_j$, ${\cal S}_j$ and $\Psi_j$,
where  ${\cal A}_1$ and ${\cal S}_1$ are $m \times m$ matrices, ${\cal A}_2$ and ${\cal S}_2$ are $\vk \times \vk$ matrices,
$\Psi_1$ is an $m\times 2$ matrix,  $\Psi_2$ is a $\vk \times 2$ matrix. Then we have
$w_{{\cal A}_2}(\la) w_{{\cal A}_1}(\la) =w_A(\la) $, where $w_A$ is the  transfer matrix function corresponding to the $S$-node
of the form
\begin{equation}       \label{t0}
A=\left[
\begin{array}{lr}
{\cal A}_1& 0\\ R & {\cal A}_2
\end{array}
\right], \quad S=\left[
\begin{array}{lr}
{\cal S}_1& 0 \\ 0 & {\cal S}_2
\end{array}
\right], \quad \Pi=\left[
\begin{array}{c}
\Psi_1  \\  \Psi_2
\end{array}
\right]; \quad R=\Psi_2\breve J\Psi_1^*{\cal S}_1   ^{-1}.
\end{equation}
\end{Pn}
\begin{proof}.
From ${\cal A}_j{\cal S}_j-{\cal S}_j{\cal A}_j^*=\Psi_j\breve J\Psi_j^*$ ($p=1,2$) and representation   (\ref{t0}),
where $R=\Psi_2\breve J\Psi_1^*{\cal S}_1   ^{-1}$, it follows that  $AS-SA^*=\Pi \breve J\Pi^*$.
Now, from   (\ref{2.22}),  (\ref{2.23})  and  (\ref{t0})  we have
$w_A(\la) =w_{2}(\la) w_{{\cal A}_1}(\la) $. 
By   (\ref{2.21}),  (\ref{2.23}),  (\ref{2.26}) and  (\ref{t0})   we obtain $T_{22}^{-1}={\cal S}_2$,
$\pi_2=\Psi_2$ and, finally, $w_2=w_{{\cal A}_2}$.    
\end{proof}
Assume for simplicity that $q_0$ in  (\ref{1.1}) is bounded in the neighbourhood of zero.
Using Corollary \ref{Pn1} and Proposition \ref{fact} we shall prove in the next section
the following theorem.
\begin{Tm}  \label{Tm1}
Let the initial system (\ref{1.1}), where the coefficient $q_0(x)$  is bounded
in the neighbourhood of zero,   be given.
Let the parameter matrices $A$, $S(0)$, and $\Pi(0)$ have the block form
\begin{equation}       \label{t1}
A=\left[
\begin{array}{lr}
{\cal A}_1& 0 \\ R & {\cal A}_2
\end{array}
\right], \quad S(0)=\left[
\begin{array}{lr}
{\cal S}_1(0) & 0 \\ 0 & 0
\end{array}
\right], \quad \Pi(0)=\left[
\begin{array}{c}
\Psi_1(0)  \\  \Psi_2(0)
\end{array}
\right], 
\end{equation}
where   ${\cal A}_1$ and ${\cal S}_1(0)$ are 
$m \times m$ matrices, $A$ and $S(0)$ are $(m+\vk) \times (m+ \vk)$ matrices $(\vk>0)$,
$\Psi_1(0)$ is an $m\times 2$ matrix,  $\Psi_2(0)$ is a $\vk \times 2$ matrix,
${\cal A}_2$ is a lower
triangular $\vk \times \vk$ matrix, and $R=\Psi_2(0)\breve J\Psi_1(0)^*{\cal S}_1   (0)^{-1}$.
Let the relations
\begin{equation}       \label{t1'}
{\cal A}_1{\cal S}_1(0)-{\cal S}_1(0){\cal A}_1^*=\Psi_1(0)\breve J\Psi_1(0)^*, \quad  \Psi_2(0)\breve J\Psi_2(0)^*=0,   \quad {\cal S}_1(0)>0
\end{equation}
hold. Introduce $\Pi(x)$ and $S(x)$ by  (\ref{2.2}) and  (\ref{2.2'}), respectively, and assume $S(x)>0$ for $x>0$.
Put $h: =[1 \quad 0 \quad 0 \, \ldots \, 0]\Psi_2$ and  suppose 
\[
h(0)=c[1\quad 0]   \quad {\mathrm{or}} \quad  h(0)=c[0\quad 1] \quad (c \not=0).
\]
Then the transformed system
\begin{equation}       \label{t2}
\Big(\frac{d}{dx}+\la q_1 +\wt q_0(x)\Big)\wt u(x, \lambda )=0,
\end{equation}
where $\wt q_0$ is given by (\ref{2.7}), has the fundamental solution
$\wt u=w_A u$.  Moreover, the coefficient $\wt q_0$ admits representation
\begin{equation}       \label{t3}
\wt q_0(x)=\frac{\k}{x}+\Up(x),
\end{equation}
where $\Up$ is bounded  in the neighbourhood of zero.
Here, $\k=\vk$, if $\vk$ is odd and $h(0)=c[0\quad 1]$ or  if $\vk$ is even and $h(0)=c[1\quad 0]$.
We have $k=-\vk$, if $\vk$ is even and $h(0)=c[0\quad 1]$ or  if $\vk$ is odd and $h(0)=c[1\quad 0]$.
\end{Tm}

\subsection{Proof of Theorem \ref{Tm1}} \label{proof}

When $\Pi(x)$ is squarely integrable in the neighbourhood of $0$, using (\ref{2.2'}) we get
$S(x)=S(0)+\int_0^x\Pi (t)\Pi (t)^*dt$.
When $S(0)=0$, we have
\begin{equation}       \label{2.38'}
S(x)=\int_0^x\Pi (t)\Pi (t)^*dt.
\end{equation}
\begin{La}  \label{Pn6} 
Let $n=1$, $A\in \BC$, $S(0)=0$, and $\Pi(0)=\a_1 [1 \quad \a_2]\breve K^*$, where $\a_1  \not=0$
and $\a_2=\pm 1$.
Assume that  the potential $q_0$ of the initial system (\ref{1.1})  is \\
(i) bounded
on $(0,\,\ve)$ for some $\ve>0$ \\or \\ (ii) has a bounded limit, when $x$ tends to $+0$. \\  
Then for the GBDT tranformation of $q_0$ defined by  (\ref{2.7})  we have
\begin{equation}       \label{2.39}
\wt q_0(x)=\Up_{\pm}(x) \mp \frac{1}{x}\s_3,
\end{equation}
where $\Up_+(x)$ and $\Up_-(x)$ are  bounded in the neighbourhood of zero, if (i) is fulfilled,
and have a bounded limit, when $x$ tends to $+0$, if (ii) is true.
\end{La}
\begin{proof}. Note that under conditions of  the lemma
we have
$AS(0)-S(0)A^*=0$.  According to   (\ref{2.13}) we have also $\Pi(0)\breve J\Pi(0)^*=0$.  So,    the identity (\ref{2.3}) holds at $x=0$,
and therefore it is satisfied at all $x\geq 0$.
Hence, the matrix function $\wt q_0$ given by  (\ref{2.7}) is the GBDT transformation 
of  $q_0$.

By the definition of $\breve K$ in  (\ref{2.13}) we get
\begin{equation}       \label{2.40}
\Pi(0)=c[1 \quad 0] \,\, {\mathrm{for}} \,\, \a_2=1, \quad  \Pi(0)=c[0 \quad 1]
\,\, {\mathrm{for}} \,\, \a_2=-1 \quad (c \not=0).
\end{equation}
By  (\ref{2.2}) one can see that
\begin{equation}       \label{2.40'}
\Pi(x) -\Pi(0)=xf(x),
\end{equation}
where  $f$ is bounded
on $[0,\,\ve)$, if $q_0$ is bounded, and  $f$ has a bounded limit for  $x \to +0$,
if  $q_0$ has a bounded limit.
From   (\ref{2.38'}),   (\ref{2.40}), and  (\ref{2.40'}) it follows that
\begin{equation}       \label{2.41}
S(x) =|c|^2x+x^2f_1(x),
\end{equation}
where  $f_1$ is also bounded
on $[0,\,\ve)$, if $q_0$ is bounded, and  $f_1$ has a bounded limit for  $x \to +0$,
if  $q_0$ has a bounded limit.
Taking into account (\ref{2.7}) and  (\ref{2.40})-(\ref{2.41})   we obtain
the statement of the lemma.
\end{proof}
Recall the notations from Proposition \ref{fact} and
consider the case $n>1$, $n_1=1$, $n_2=n-1$, so that $\pi_1$ is the first block (i.e. the first row) of $\Pi$,
and $\pi_2$ is given by  (\ref{2.26}).
\begin{La}  \label{La7}
Let  the initial system (\ref{1.1})  be given
and let  the  potential $q_0$ be bounded in the neighbourhood of zero.
Suppose that $A$ is a lower triangular $n \times n$ $(n > 1)$ matrix,
that relations  (\ref{2.2}) and (\ref{2.2'})  hold and that
\begin{equation}       \label{2.51}
S(0)=0, \quad \Pi(0)\breve J\Pi(0)^*=0, \quad S(x)>0 \quad {\mathrm{for}} \quad x>0.
\end{equation}
Then,  we have 
\begin{eqnarray}       \label{2.52}
&&\wh q_0(x)=\wh \Up_{+}(x) - \frac{1}{x}\s_3 \quad {\mathrm{for}} \quad \pi_1(0)=c[1 \quad 0];
\\  \label{2.53}
&& \wh q_0(x)=\wh \Up_{-}(x) + \frac{1}{x}\s_3 \quad {\mathrm{for}} \quad \pi_1(0)=c[0 \quad 1],
\end{eqnarray}
where $c\in\BC\backslash \{0\}$, $\wh q_0$ is given by   (\ref{2.29}) and $\wh \Up_{\pm}$ is bounded in the neighbourhood of zero.
Moreover, the matrix function $\pi_2(x)$ is continuous at zero and the relations
\begin{equation}       \label{2.54}
\lim_{x\to +0}\big(T_{22}^{-1}\big)(x)=0, \quad \big(T_{22}^{-1}\big)(x)>0  \quad {\mathrm{for}} \quad x>0,
 \quad \pi_2(0)\breve J\pi_2(0)^*=0
\end{equation}
hold for both cases   (\ref{2.52}) and  (\ref{2.53}). The matrix functions $T_{22}^{-1}$
and $\pi_2$ satisfy  (\ref{2.28}).
\end{La}
\begin{proof}. 
First, compare  Proposition \ref{fact} and Lemma \ref{Pn6}  to see that $\wh q_0$ from
(\ref{2.29}) coincides with
$\wt q_0$ from  (\ref{2.39}). Thus,  (\ref{2.52}) and  (\ref{2.53}) are immedate from
Lemma  \ref{Pn6}. 

As $S(x)>0$, so we get $T(x)>0$,  $T_{22}(x)>0$ and, finally,  $\big(T_{22}^{-1}\big)(x)>0$.
Similar to the case treated  in
Lemma \ref{Pn6}
 we have
$AS(0)-S(0)A^*=0$ and $\Pi(0)\breve J\Pi(0)^*=0$, i.e.,   (\ref{2.3}) holds at $x=0$.
Hence, the conditions of  Proposition \ref{fact}  hold. By 
Proposition \ref{fact} $T_{22}^{-1}$
and $\pi_2$ satisfy  (\ref{2.28}). 

Recall that $S(0)=0$. Therefore, using formulae  (\ref{2.38'}) and (\ref{2.41})
we derive
\begin{equation}       \label{2.55}
S_{22}(0)=0, \quad  \lim_{x\to +0}S_{21}(x)S_{11}(x)^{-1}S_{12}(x)=0.
\end{equation}
The second equality in  (\ref{2.21}) and formula (\ref{2.55}) imply the first relation in 
(\ref{2.54}). From the definition of $\pi_2$ in   (\ref{2.26})
 and from the first relation
in  (\ref{2.21}) it follows that 
\begin{equation}       \label{2.55d}
\pi_2=[-S_{21}S_{11}^{-1} \quad I_{n-1}]\Pi.
\end{equation}
Now, use again (\ref{2.38'}) and (\ref{2.41}) to see that $\pi_2$ has a limit,
when $x$ tends to $+0$. Therefore, the third relation in (\ref{2.54}) 
is immediate from the first relations in  (\ref{2.26}) and (\ref{2.54}).
\end{proof}
 Next, consider the case of the initial system with $\k\not=0$,
where
$V \in L^p_{2 \times 2}(0, \ve)$,  that is, the entries of $V$ belong
to $L^p$ in the neighbourhood of zero.

\begin{La}  \label{La8}
Let the initial Dirac system  (\ref{1.1}) be given, where
\begin{equation}       \label{2.55'}
q_0(x)=i\s_2\Big(V(x)+\frac{\k}{x}\s_1\Big)=i\s_2V(x)+\frac{\k}{x}\s_3, \quad \k\not=0,
\end{equation}
 and $V \in L^p_{2 \times 2}(0, \ve)$.
Suppose that $A$ is a lower triangular $n \times n$ $(n > 1)$ matrix,
that relations  (\ref{2.2}) and (\ref{2.2'})  hold, that $S$ and $\Pi$ are continuous at zero and that
 (\ref{2.51}) is true. Then, putting $n_1=1$, we get
\begin{equation}       \label{2.56}
\wh q_0(x)=-\frac{1+\k}{x}\s_3+\wh \Up_+(x) \, \,{\mathrm{for}} \, \k>0, \quad
\wh q_0(x)=\frac{1-\k}{x}\s_3+\wh \Up_-(x) \, \, {\mathrm{for}} \, \k<0,
\end{equation}
where  $\wh q_0$ is given by   (\ref{2.29}) and  $\wh \Up_{\pm}\in L^p_{2 \times 2}$    in the neighbourhood of zero.
Moreover, the matrix function $\pi_2(x)$ is continuous at zero and the relations  (\ref{2.54}) 
hold. The matrix functions $T_{22}^{-1}$
and $\pi_2$ satisfy  (\ref{2.28}).
\end{La}
To prove this  lemma we shall need  Lemma A1 from  \cite{AHM}:
\begin{La}  \label{La9} Let $q_0$ be given by  (\ref{2.55'}).
Assume $\k \in \BZ \backslash{0}$, $\la \in \BC$, and  $V \in L^p_{2 \times 2}$   $(1 \leq p <\infty)$
in the neighbourhood of zero. Then, there are two types of   nontrivial solutions
$y={\mathrm{col}}[y_1 \quad y_2]$
of  (\ref{1.1}) . Namely,
as $x \to +0$, either the limit
$\lim_{x\to +0}x^{-\k}y(x,\la)\not=0$ exists and $x^{-\k-1}y_1(x,\la) \in L^p$ in the neighbourhood of zero
or  the limit
$\lim_{x\to +0}x^{\k}y(x,\la)\not=0$ exists and $x^{\k-1}y_2(x,\la) \in L^p$ in the neighbourhood of zero.
\end{La}
A similar result is true for the bounded $V$.
\begin{La}\label{LaInf}
Assume that $\k \in \BZ \backslash{0}$, $\la \in \BC$, and $V\in L^{\infty}_{2 \times 2}$, that is, $V$ is bounded
in the neighbourhood of zero. Then, the entries of the bounded in  the neighbourhood of zero
solution $y$ of  (\ref{1.1}) have the property:

If $\k>0$, then $x^{-\k-1}y_1(x,\la) \in L^\infty$ in the neighbourhood of zero,
and if $\k<0$, then $x^{\k-1}y_2(x,\la) \in L^\infty$ in the neighbourhood of zero.
\end{La}
\begin{proof}.
Lemma \ref{LaInf} easily follows from the proof of  Lemma A1 \cite{AHM}.
Consider, for instance, the case $\k>0$. When $V=0$, equation
(\ref{1.1}) has a solution $\phi_0={\mathrm{col}}[\phi_1^0 \quad \phi_2^0]$, such that
the functions $x^{-\k-1}\phi_1^0(x,\la) $ and  $x^{-\k}\phi_2^0(x,\la) $ are bounded
in the neighbourhood of zero,
and (\ref{1.1}) has another solution $\psi_0={\mathrm{col}}[\psi_1^0 \quad \psi_2^0]$,
such that there is a limit $\lim_{x\to +0}x^{\k}\psi_0(x,\la) \not=0$. 
By the considerations of \cite{AHM}, the bounded solution $y$ of the Dirac equation (\ref{1.1})
with a bounded $V$ satisfies the integral equation
\begin{equation}       \label{2.56d}
\big(I-{\cal{L}}\big)y=\g,  \quad  
{\cal{L}}:=\psi_0(x,\la) \int_0^x[\phi_1^0(t,\la)  \quad \phi_2^0(t,\la) ]V(t)\, \cdot \, dt, 
\end{equation}
\begin{equation}       \label{2.56d1}
 \g(x,\la) :=\phi_0(x,\la) \Big(1-\int_0^x[\psi_1^0(t,\la)  \quad \psi_2^0(t,\la) ]V(t)u(t,\la) dt\Big), \quad 0<x<\ve .
\end{equation}
It follows from the properties of $\phi_0$, Lemma \ref{La9} and definition (\ref{2.56d1}) that
\begin{equation}       \label{2.56d2}
x^{-\k-1}\g_1(x,\la) \in L^{\infty}, \quad  x^{-\k}\g_2(x,\la) \in L^{\infty} \quad(\g=[\g_1 \quad \g_2], \, 0<x<\ve).
\end{equation}
Using  (\ref{2.56d}) and (\ref{2.56d2}), by induction we get for some $C>0$ and $\ve(C)>0$ that
\begin{equation}       \label{2.56d3}
\| \Big({\cal{L}}^j\g\Big)(x)\|\leq C^{j+1}x^{j+\k}, \quad   j>0, \quad 0<x<\ve ,
\end{equation}
where $\|\cdot\|$ is the usual vector norm.
As  $y=\sum_{j=0}^{\infty}{\cal{L}}^j \g$, it follows from (\ref{2.56d2}) and (\ref{2.56d3})
that $x^{-\k-1}y_1(x,\la) \in L^\infty$.
For the case $\k<0$
the proof is similar.
\end{proof}

Notice that from 
$S(x)=\int_0^x\Pi(t)\Pi(t)^*dt>0$
one easily gets $\pi_1(x)\not\equiv 0$ in any neighbourhood  of zero.
By (\ref{2.2}) and by the third equality in  (\ref{r16})   we see
that $\breve J^*\pi_1^*$ satisfies  (\ref{1.1}), where $\la=A_{11}^*$.
Therefore,
the next corollary immediately follows from Lemmas \ref{La9} and \ref{LaInf}.
\begin{Cy}  \label{Cy10}
Let the conditions of Lemma   \ref{La8} be fulfilled, where $1\leq p \leq \infty$.
Then for the entries $h_1$ and $h_2$ of $\pi_1=[h_1 \quad h_2]$,  for some $\ve>0$
and for some $c\in \BC\backslash 0$ we have
\begin{eqnarray}       \label{2.57}
&&\lim_{x\to +0}x^{-\k}h_1(x)=c\not=0, \quad x^{-\k-1}h_2(x) \in L^p(0,\ve)  \quad {\mathrm{if }} \quad \k>0,
\\   \label{2.58}
&&\lim_{x\to +0}x^{\k}h_2(x)=c\not=0, \quad x^{\k-1}h_1(x)\in L^p(0,\ve),  \quad {\mathrm {if}} \quad \k<0.
\end{eqnarray}
\end{Cy}
\begin{proof} of Lemma  \ref{La8}. 
In a  similar to the proof of Lemma \ref{La7} way  one shows that  the conditions of Proposition \ref{fact}
are fulfilled and so (\ref{2.28}) is true.

From $S(x)>0$ it follows that  $\big(T_{22}^{-1}\big)(x)>0$ for $x>0$,
i.e., the second relation in (\ref{2.54}) is valid.
Taking into account the equality $S(0)=0$, formula (\ref{2.27}) and Corollary \ref{Cy10} we see that
\begin{equation}       \label{2.59}
S_{11}(x)=\int_0^x\pi_1(t)\pi_1(t)^*dt=(2|k|+1)^{-1}|c|^2x^{2|k|+1}+o\Big(x^{2|k|+1}\Big), \quad x\to +0.
\end{equation}
Moreover, in view of  (\ref{2.38'}) and  Corollary \ref{Cy10} we obtain
\begin{equation}       \label{2.59d}
S_{21}(x)=\int_0^xP_2\Pi(t)\pi_1(t)^*dt=(|k|+1)^{-1}x^{|k|+1}P_2\Pi(0)\pi_1(0)^*+o\Big(x^{|k|+1}\Big).
\end{equation}
Using the second relation in (\ref{2.21})  and formulae (\ref{2.59}) and (\ref{2.59d}),
we get the first relation in (\ref{2.54}). From the definition of $\pi_1$  in (\ref{2.25})
and formula (\ref{2.55d}) it follows:
\begin{equation}       \label{2.60}
\pi_2=-S_{21}S_{11}^{-1}\pi_1+P_2\Pi.
\end{equation}
By  (\ref{2.57})-(\ref{2.60}) the matrix function $\pi_2$ is continuous at zero.
Hence, the third relation in  (\ref{2.54}) follows from the first relations in (\ref{2.26})
and (\ref{2.54}).

It remains to prove (\ref{2.56}). For this purpose we shall follow a nice scheme
from \cite{AHM}.
According to  (\ref{2.57}) and  (\ref{2.59})
 we have
\begin{equation}       \label{2.61}
\breve J\pi_1^*S_{11}^{-1}\pi_1\breve J^*-\pi_1^*S_{11}^{-1}\pi_1=-\frac{h_1\ov h_1}{S_{11}}\s_3+\Up_1=-\frac{\pi_1 \pi_1^*}{S_{11}}\s_3+\Up_2
\, {\mathrm{for}} \, \k>0,
\end{equation}
where  $\Up_1, \, \Up_2 \, \in L^p_{2 \times 2}$. (In the proof of this lemma 
functions $\Up_j$ ($1\leq j \leq 8$) are considered in the neighbourhood of zero.)
 According to  (\ref{2.58}) and  (\ref{2.59})
 we have
\begin{equation}       \label{2.62}
\breve J\pi_1^*S_{11}^{-1}\pi_1\breve J^*-\pi_1^*S_{11}^{-1}\pi_1=\frac{h_2\ov h_2}{S_{11}}\s_3+\Up_3=\frac{\pi_1 \pi_1^*}{S_{11}}\s_3+\Up_4
\, {\mathrm{for}} \, \k<0,
\end{equation}
where  $\Up_j \in L^p_{2 \times 2}$  ($j=3,4$).
In view of  (\ref{2.27}),  rewrite ${\pi_1 \pi_1^*}{S_{11}^{-1}}$ in the form
\begin{equation}       \label{2.63}
\frac{\pi_1(x) \pi_1(x)^*}{S_{11}(x)}=\frac{d}{dx}\ln S_{11}(x)
=\frac{2|\k|+1}{x}+\frac{g_x(x)}{g(x)}, \quad  g(x):=x^{-2|\k|-1} S_{11}(x).
\end{equation}
From  (\ref{2.29}), (\ref{2.55'}) and (\ref{2.61})--(\ref{2.63}) it follows that
\begin{equation}       \label{2.63d}
\wh q_0(x)=-\frac{1+\k}{x}\s_3-\frac{g_x(x)}{g(x)}\s_3+\Up_5(x)  \quad {\mathrm{for}} \quad  \k>0, 
\end{equation}
\begin{equation}       \label{2.63'}
\wh q_0(x)=\frac{1-\k}{x}\s_3+\frac{g_x(x)}{g(x)}\s_3+ \Up_6(x)  \quad {\mathrm{for}} \quad \k<0,
\end{equation}
where  $\Up_j \in L^p_{2 \times 2}$  ($j=5,6$). Put
\begin{equation}       \label{2.64}
f(x):=x^{-2|k|}\pi_1(x) \pi_1(x)^*=x^{-2|k|}\pi_1(x) \breve J\breve J^*\pi_1(x)^*.
\end{equation}
Recall that $\breve J\pi^*$ satisfies Dirac equation (\ref{1.1}) with $q_0$ of the form
 (\ref{2.55'}), where the entries of $V$ belong $L^p$. Hence, differentiating
 (\ref{2.64}) and taking into account Corollary \ref{Cy10} we get
\begin{equation}       \label{2.65}
\frac{d}{dx}f(x):=-\frac{2|\k|}{x^{2|k|+1}}\pi_1(x) \pi_1(x)^*+x^{-2|k|}\pi_1(x) \breve J\Big(-\frac{2\k}{x}\s_3\Big)\breve J^*\pi_1(x)^*+\Up_7(x),
\end{equation}
where $ \Up_7\in L^p$. Apply again Corollary \ref{Cy10} to the right-hand side of (\ref{2.65}) to obtain
\begin{equation}       \label{2.66}
\frac{d}{dx}f(x)=-\frac{2|\k|}{x^{2|k|+1}}\pi_1(x) \pi_1(x)^*+\frac{2\k}{x^{2|k|+1}}\pi_1(x) \s_3\pi_1(x)^*+\Up_7(x)=\Up_8(x),
\end{equation}
where $\Up_8\in L^p$.  Next, use (\ref{2.59}),  (\ref{2.64}) and (\ref{2.66})  to   rewrite $S_{11}$ in the form
\begin{equation}       \label{2.67}
S_{11}(x)=\int_0^x t^{2|\k|}f(t)dt=\frac{x^{2|k|+1}}{2|k|+1}f(x)- \frac{1}{2|k|+1}\int_0^x  t^{2|\k|+1}\frac{df}{dt}(t)dt.
\end{equation}
Taking into account (\ref{2.63}) and (\ref{2.67}), we have
\[
g(x)= \frac{1}{2|k|+1}\Big(f(x)-    {x^{-2|k|-1}} \int_0^x  t^{2|\k|+1}\frac{df}{dt}(t)dt\Big), 
\]
which implies
\begin{equation}       \label{2.69}
\frac{d}{dx}g(x)=\big({\cal K}f_x\big)(x), \quad {\cal K}:={x^{-2|k|-2}} \int_0^x  t^{2|\k|+1}\, \cdot\,dt.
\end{equation}
As ${\cal K}$ is a bounded in $L^p(0,\ve)$ ($1\leq p\leq \infty$) Hardy operator, formulae   (\ref{2.66}) and (\ref{2.69})
imply
\begin{equation}       \label{2.70}
\frac{dg}{dx}(x) \in L^p(0, \ve).
\end{equation}
From  (\ref{2.63d}), (\ref{2.63'}), and (\ref{2.70}) follow relations (\ref{2.56}).
\end{proof}
\begin{proof} of Theorem \ref{Tm1}.
Using  (\ref{t1}) and (\ref{t1'})   one can see
that  (\ref{2.4}) holds at $x_0=0$.
As $q_0$ is bounded, the condition of Remark \ref{Rk1}
is satisfied and identity (\ref{2.4})  at $x_0=0$ implies (\ref{2.3}).
Hence, the conditions of Proposition \ref{fact}
are fulfilled. Therefore,
the GBDT generated by the parameter matrices
 $A$, $S(0)$, and $\Pi(0)$   is a superposition
of two GBDTs, where the first GBDT is  generated by ${\cal A}_1    $, ${\cal S}_1   (0)$, and $\Psi_1(0)$
(in addition to Proposition \ref{fact} see also (\ref{2.22})--(\ref{2.23'}) and Remark \ref{sup}).
As ${\cal S}_1(0)>0$ and $q_0$ is bounded
in the neighborhood  of zero,
so
the transformation of $q_0$, which is generated by ${\cal A}_1    $, ${\cal S}_1   (0)$, $\Psi_1(0)$,  is bounded too. 
Denote this transformation
by $q_{0}^{(1)}$.

The second GBDT in the superposition is determined by the $\vk \times \vk$ matrix ${\cal A}_2$ 
and by the matrix functions
$T_{22}^{-1}(x)$, and $\pi_2(x)$.
 As $S_{21}(0)=0$ and ${\cal S}_1(0)>0$
formula (\ref{2.55d}) implies $\lim_{x\to +0}\pi_2(x)=\Psi_2(0)$.
According to the second relation in  (\ref{2.21})
we have $\lim_{x\to +0}T_{22}^{-1}(x)=0$, that is, we may put $T_{22}^{-1}(0)=0$.
Thus, the second GBDT is generated by the parameter matrices
${\cal A}_2   $, $T_{22}^{-1}(0)=0$, and $\pi_2(0)=\Psi_2(0)$, which satisfy
the identity ${\cal A}_2T_{22}^{-1}(0) -T_{22}^{-1}(0)  {\cal A}_2^*= \Psi_2(0)\breve J\Psi_2(0)^*$.

If $\vk=1$, then the   potential  $q_{0}^{(1)}$ and the $S$-node ${\cal A}_2   $, $T_{22}^{-1}(0)=0$, and $\pi_2(0)=\Psi_2(0)$
satisfy  the
conditions of Lemma \ref{Pn6} and the statement of the theorem is true.

If  $\vk>1$, then 
${\cal A}_2   $, $T_{22}^{-1}(0)=0$ and $\pi_2(0)$ satisfy  the
conditions of Lemma \ref{La7}, and the second GBDT is itself a superposition of GBDTs.
Taking into account the block
representation
\begin{equation}       \label{2.71}
{\cal A}_2= \left[
\begin{array}{lr}
{\cal A}^{(1)}&0 \\  \, \,* &{\cal A}_2^{(2)}
\end{array}
\right], \quad \pi_2= \left[
\begin{array}{c}
\pi^{(1)}\\ *
\end{array}
\right], \quad {\cal A}^{(1)}, \, \pi^{(1)} \, \in \BC
\end{equation}
we derive that the transformation of $q_{0}^{(1)}$ generated by ${\cal A}_2   $, $0$, and $\pi_2(0)$
is a superposition of the transformation generated by ${\cal A}^{(1)}$, $0$, and $ \pi^{(1)}(0)$, 
and  of the transformation generated by ${\cal A}_2^{(2)}$, $0$, and $\pi_2^{(2)}(0)$, respectively, where $\pi_2^{(2)}$
is constructed analogously to the construction of $\pi_2$ in Proposition \ref{fact}.
 (The transformation of $q_{0}^{(1)}$  generated by
 ${\cal A}^{(1)}$, $0$, and $ \pi^{(1)}(0)$ is denoted by $q_{0}^{(2)}$.)
Moreover, for $\vk >2$ the
potential $q_{0}^{(2)}$ and the $S$-node 
${\cal A}_2^{(2)}$, $0$, and $\pi_2^{(2)}(0)$ satisfy the conditions  of Lemma \ref{La8}, 
and according to Lemma \ref{La8} this holds for further $\vk -3$ steps too.
At each step  we obtain a new coefficient $q_{0}^{(s)}$. The absolute value of $\k^{(s)}$ 
in the representation (\ref{2.55'}) of $q_{0}^{(s)}$ ($s>1$) equals $s-1$, and the sign of $\k$ changes at each step
till we come to the final transformation.

If $\vk>1$, the final transformation is generated by
\[
{\cal A}_2^{(\vk)}=({\cal A}_2)_{\vk\vk}\in \BC, \quad 0, \quad \pi_2^{(\vk)}(0)=(0,0).
\]
Though we do not apply Lemma  \ref{La8} at this last step directly, the
proof of the representation of $\wh q_0$ from Lemma  \ref{La8}  remains valid for $\wt q_0= q_0^{(\vk+1)}$.
\end{proof}

\section{Direct and inverse problems}\label{dirinv}
\setcounter{equation}{0}
In this section we  consider transformed systems constructed in Section \ref{Prel},
assuming that the initial systems are trivial, that is, either  $v \equiv 0$ or $\xi \equiv 0$.
We use the notion of a minimal  relization
and some other notions from system theory, which are defined
in Appendix A.
\subsection{Dirac systems}
\subsubsection{Self-adjoint Dirac system}
First, consider the self-adjoint Dirac system constructed in Proposition \ref{PnGTDS}
\begin{equation}       \label{i1}
\frac{d}{dx}\wt u(x, \la)=i\big( \la j+j \wt V(x)\big)\wt u(x,\la) \quad (0 \leq x <\infty),
\quad \wt V= \left[
\begin{array}{cc} 0 &
\wt v \\ \wt v^* & 0
\end{array}
\right].
\end{equation}
We assume that
\begin{equation}       \label{i2}
 S(0)=I_n>0, \quad A-A^*=i\Pi(0)j\Pi(0)^*, \quad v(x) \equiv 0,
\end{equation}
where $v$ is the potential of the initial Dirac system (\ref{0x1.-4}).
Partition $\Pi$ into two $n \times p$ blocks. Using (\ref{r5}), where
$V \equiv 0$, we derive
\begin{equation}       \label{i3}
\Pi(x)=[\Phi_1(x) \quad \Phi_2(x)]=[e^{-ixA}\Phi_1(0) \quad e^{ixA}\Phi_2(0)].
\end{equation}
It follows from  (\ref{se.53}), (\ref{se.65}), and (\ref{i2}) that
\begin{eqnarray}       \label{i4}
&&\wt v(x)=-2i  \Phi_1(0)^*e^{ixA^*}S(x)^{-1}e^{ixA}\Phi_2(0), \\
&&      \label{i5}
S(x)=I_n+\int_0^x \Pi(t)\Pi(t)^*dt>0.
\end{eqnarray}
As $S(0)=I_n$ is fixed, parameter matrices in (\ref{i3})--(\ref{i5}) are $A$ and $\Pi(0)$
or, equivalently, $A$, $\Phi_1(0)$, and  $\Phi_2(0)$.
\begin{Dn}\label{PE} \cite{GKS1}
The potentials $\wt v$ of the form (\ref{i4}) are called pseudo-exponential
and the class of such potentials is denoted by PE.
\end{Dn}
If  $\wt V$ is locally summable on $[0, \infty)$, there is a unique
Weyl function of system (\ref{i1})  (see, for instance, \cite{SaA7}),
which is defined in the following way.
\begin{Dn} \label{DnW1}
A holomorphic $p \times p$ matrix function $\varphi$ 
such that
\begin{eqnarray} \label{i6}
&& \int_0^\infty \left[ \begin{array}{lr} I_p &  i \varphi (\la) ^*
\end{array} \right]
  K \wt u(x, \la) ^*
\wt  u(x, \la) K^*
 \left[ \begin{array}{c}
I_p \\ - i \varphi (\la)  \end{array} \right] dx < \infty , 
\\ && \label{i7}
\la\in {\BC}_+, \quad \wt u(0,\la)=I_m, \quad K:=   \frac{1}{\sqrt{2}}       \left[
\begin{array}{cc} I_p &
-I_{p} \\ I_{p} & I_p
\end{array}
\right]
\end{eqnarray}
is called a Weyl function  of Dirac-type system (\ref{i1}) on
$[0, \, \infty)$.
\end{Dn}
Moreover, the function $\vp$, which satisfies (\ref{i6}), is unique
even without the analyticity requirement.

Now, we  shall again consider the case $v \in$ PE.
As $V=0$, we have $u(x,\la)=\exp(ix\la j)$ in (\ref{r19}), and so formula (\ref{r19}) takes the form
\begin{equation}       \label{i8}
\wt u(x,\la)=w_A(x, \la)e^{ix \la j}w_A(0, \la)^{-1}.
\end{equation}
By Proposition  6.2 \cite{GKS6}  for every $\la \in \BC$, excluding a finite set of points,
there is a limit
\begin{equation}       \label{i9}
\lim_{x \to \infty}w_A(x, \la)=f_A(\la).
\end{equation}
If the equality
\begin{equation}       \label{i10}
K^*
 \left[ \begin{array}{c}
I_p \\ - i \varphi (\la)  \end{array} \right] =w_A(0, \la)\left[ \begin{array}{c}
I_p \\ 0 \end{array} \right]c(\la)
\end{equation}
holds for some $p \times p$  matrices $c(\la)$ and $\vp(\la)$ ($\la \in \BC_+$), then
formulas  (\ref{i8}) and (\ref{i9}), and exponential decay of  $e^{ix\la}$ imply
inequality (\ref{i6}). As $K^*=K^{-1}$, it is easy to see that 
(\ref{i10}) is everywhere, excluding a finite set, equivalent to the relation
\begin{equation}       \label{i11}
\vp(\la)=i\left([0 \quad I_p]Kw_A(0, \la)\left[ \begin{array}{c}
I_p \\ 0 \end{array} \right]\right)\left([I_p \quad 0]Kw_A(0, \la)\left[ \begin{array}{c}
I_p \\ 0 \end{array} \right]\right)^{-1}
\end{equation}
Taking into account definitions (\ref{se.56}) and (\ref{i7}) of $w_A$ and $K$,
respectively, and formula (\ref{app4}), one  gets the following result (see \cite{GKS1} and Theorem 5.1 in 
\cite{GKS6}).
\begin{Tm} \label{DW}
 Let $\wt v \in {\mathrm{PE}}$. Then the Weyl
function of the Dirac system (\ref{i1}) admits realization
\begin{equation}       \label{i12}
\vp(\la)=iI_p+2\Phi_2(0)^*(\la I_n -\breve A)^{-1}\Phi_1(0), \quad \breve A:=A-i\Phi_1(0)
\big(\Phi_1(0)+\Phi_2(0)\big)^*.
\end{equation}
\end{Tm}
Note that the  Weyl function $\vp$ is a Herglotz function and admits representation
\begin{equation}       \label{i13}
\vp(\la)=\nu+\int_{-\infty}^{\infty}\Big(\frac{1}{t-\la}-\frac{t}{1+t^2}\Big)d \tau(t), \quad (\nu=\nu^*, \,\, \tau \uparrow).
\end{equation}
Here $\tau$ is the spectral function of system  (\ref{i1}). See Theorem 4.3 \cite{GKS1}, where this fact
is proved and an explicit expression for $\tau$ in terms of the parameter matrices is given.
The inverse problem to recover $\wt v \in $ PE from $\tau$   is solved explicitly in Theorem 4.5 \cite{GKS1}.
(The corresponding bound states are constructed in \cite{GKS4}.)
The inverse problem to recover $\wt v$ from the spectral density for the particular case $\s(A) \subset \BC_+$
was treated earlier in \cite{AG1}. 

The left reflection coefficient for system (\ref{i1}), where $\wt v \in$ PE, is expressed
via  the Weyl function by the formula (see p. 33 in \cite{GKS6}):
\begin{equation}       \label{i14}
R_L(\lambda) =  -\big(I_p + i\varphi(\lambda)\big) \big(I_p
-i\varphi(\lambda)\big)^{-1}.
\end{equation}
Recall that $\vp$ has properties
\begin{equation}       \label{i15}
\Im
\vp(\la) \geq 0 \quad (\la \in \BC_+), \quad \lim_{\la \to \infty}\vp(\la)=iI_p.
\end{equation}
It follows from (\ref{i14}) and (\ref{i15}) that
\begin{equation}       \label{i16}
\|
R_L(\la)\| \leq 1  \quad (\la \in \BC_+), \quad \lim_{\la \to \infty}R_L(\la)=0.
\end{equation}

To recover system  (\ref{i1}) from its Weyl function we should recover
parameter matrices $A$ and $\Phi_k(0)$ ($k=1,2$). First, we recover
$R_L$ via formula (\ref{i14}).
By the second relation in (\ref{i16}) $R_L$ is a strictly proper matrix function
and admits a minimal realization 
\begin{equation}       \label{i17}
R_L(\la)={\cal C}(\lambda I_n-{\cal A})^{-1}{\cal B},
\end{equation}
where ${\cal C}$ is a $p \times n$ matrix and  ${\cal B}$ is an $n \times p$ matrix.
According to the first relation in (\ref{i16}), $R_L$  admiting the minimal realization (\ref{i17})
is also contractive. Hence, by Theorems 21.1.1 and 21.1.3 in \cite{LR}
the Riccati equation
\begin{equation}       \label{i18}
{\cal X}{\cal C}^*{\cal C}{\cal X}-i(   {\cal A}{\cal X} - {\cal X}{\cal A}^*)+ {\cal B}{\cal B^*}=0
\end{equation}
has a positive solution ${\cal X}>0$.
Put
\begin{equation}       \label{i19}
A={\cal X}^{-\frac{1}{2}}{\cal A}{\cal X}^{\frac{1}{2}}+
i{\cal X}^{-\frac{1}{2}}{\cal B}{\cal B}^*{\cal X}^{-\frac{1}{2}}, \quad
\Phi_1(0)={\cal X}^{-\frac{1}{2}}{\cal B},  \quad \Phi_2(0)=-i{\cal X}^{\frac{1}{2}}{\cal C}^*.
\end{equation}

The second relation in (\ref{i2}) is immediate  from (\ref{i18}) and (\ref{i19}),
that is, $A$ and $\Phi_k(0)$ ($k=1,2$) are some parameter matrices
and we can apply Theorem \ref{DW}. It easily follows
(see Section 9 in \cite{GKS6}) that $\vp$ given by  (\ref{i12}) and (\ref{i19})
coincides with the  $\vp$ from which $R_L$ was recovered.
\begin{Tm} \label{IW} \cite{GKS6}
Let $\vp$ be a rational $p \times p$ matrix function, which satisfies (\ref{i15}).
Then $\vp$ is the Weyl function of some system  (\ref{i1}), where
$\wt v \in {\mathrm{PE}}$. To recover $\wt v$ we take a minimal realization
(\ref{i17})
of $R_L$ given by (\ref{i14}) and define parameter matrices by (\ref{i19}),
where ${\cal X}>0$ is a solution of the Riccati
equation (\ref{i18}). After that we apply formulas (\ref{i3})--(\ref{i5}).
\end{Tm}
By  Theorem 5.4 \cite{GKS6} all the potentials $\wt v \in$ PE admit representation
(\ref{i4}), where parameter matrices have additional properties
\begin{equation}       \label{i20}
\s(A)\subset \ov \BC_+, \quad
 {\mathrm{span}}\bigcup_{k=0}^{n-1}\im
{A}^k\Phi_1(0)=\BC^n,\quad {\mathrm{span}} \bigcup_{k=0}^{n-1}\im  A^k\Phi_2(0)=\BC^n,
\end{equation}
and  Im means image. Hence, without loss of generality we can assume (\ref{i20}).

\subsubsection{Potentials with singularities}
Now, consider a more general case of Dirac systems,
including systems with singularities. Namely, let
\begin{equation}  \label{i21}
j = \left[
\begin{array}{cc}
I_{p_1} & 0 \\ 0 & -I_{p_2}
\end{array}
\right], \quad 
S(0)=S_0, \quad A S_0-S_0A^*=i\Pi(0)j\Pi_0^*,
\end{equation}
where $p_1, \, p_2>0$, $p_1+p_2=m$. Formula
(\ref{i5}) takes the form
\begin{equation}       \label{i22}
S(x)=S_0+\int_0^x \Pi(t)\Pi(t)^*dt>0.
\end{equation}
As before, $\Pi$, $w_A$, and $\wt v$ are defined via formulas (\ref{i3}), (\ref{se.56}), and (\ref{i4}),
respectively, and we assume that (\ref{i20}) holds.  Here $\Phi_1$, $\Phi_2$, and $\wt v$
are $n \times p_1$,  $n \times p_2$, and $p_1 \times p_2$ matrix functions,
respectively.
\begin{Dn}\label{GPE} \cite{FKS}
The potentials $\wt v$ of the form (\ref{i4}),
where (\ref{i20}) and (\ref{i21}) hold, and  $S$ is given by (\ref{i22}), are called 
generalized pseudo-exponential
and the class of such potentials is denoted by GPE.
\end{Dn}

Dirac system with $\wt v \in$ GPE was treated on the whole axis
in \cite{SaA8}, and results on supertransparent potentials and soliton-positon
interactions were obtained. In particular, it was shown in the proof
of  Theorem 2.1 from \cite{SaA8} that for some $x_0 \in \BR$
\begin{equation}       \label{i23}
S(x)>0 \quad {\mathrm{for}} \quad x >x_0.
\end{equation}
If $S_0 \not> 0$, then $\det \, S(x_k)=0$ in some point (or points) on $\BR_+\cup\{0\}$,
which means that the GPE class includes potentials with singularities.
Nevertheless, as $S(x)$ admits continuation meromorphic in $x$, the determinant
$\det \, S(x)$ turns into zero only in a finite number of points on $\BR_+\cup\{0\}$.
The proof that the matrix function
\begin{equation}       \label{i24}
\wt u(x,\la)=w_A(x,\la)e^{ix\la j}
\end{equation}
satisfies (\ref{i1}) and is nondegenerate (excluding a finite number of 
values of $\la$ and zeros of $S(x)$)
remains true for Dirac system with $j$ given by (\ref{i21}).
Moreover, $\wt u$ as well as $S$ is meromorphic in $x$. Therefore,
we call $\wt u$ a fundamental solution of  (\ref{i1}),
which agrees with the standard \cite{M2} requirement
for the fundamental and special solutions in the case of singularities
to be defined by the same formula on the whole domain.

The asymptotics of $w_A$ in (\ref{i24}) under condition (\ref{i20}) is given in the next proposition
(see Proposition 3.1 in \cite{FKS} or Theorem 3.1 in \cite{SaA8}).
\begin{Pn} \label{PnAs} Let
$A$, $S(0)$, and $\Phi_k(0)$ $(k=1,2)$
satisfy (\ref{i20}) and (\ref{i21}). Then there is a limit
\begin{equation} \label{i24'}
 \om= \lim_{x \to \infty}\Big(e^{-ix A } S(x)e^{ix A^{*} }\Big)^{-1} \geq 0,
\end{equation}
 and we have
\begin{equation}\label{i24''}
\lim_{x \to \infty}w_{A}(x,\lambda )=
\left[\begin{array}{lr} I_{p_1}
 & 0 \\
0 &   \chi(\la) \end{array}\right], \quad \chi(\la):=
I_{p_2}+i\Phi_2(0)^*\om(A-\la I_n )^{-1}\Phi_2(0).
\end{equation}
\end{Pn}
Transmission and reflection coefficients for system (\ref{i1})
are defined \cite{FKS, GKS6, YL} in terms of   special solutions of  (\ref{i1}), that is,
$m \times p_1$ and $m
\times p_2$ solutions ${\cal Y}$ and ${\cal Z}$, respectively, 
which  we determine by the relations
\begin{equation}\label{i25}
 {\cal Y}(x,\la )=e^{ix\la}\left[ \begin{array}{c} I_{p_1} \\ 0
\end{array}
  \right] +o(1)
\quad (x\to\infty ),\quad {\cal Z}(0,\la )=\left[
\begin{array}{c} 0 \\
I_{p_2} \end{array} \right], \quad \la \in \BR.
\end{equation}
  Put
\begin{eqnarray}
&&T_L(\la ):={\cal Y}_1(0,\la )^{-1},\quad R_L(\la):={\cal Y}_2(0,\la
)\, {\cal Y}_1(0,\la
)^{-1},\label{i26}\\ \noalign{\vskip6pt} && R_R(\la ):=
\Bigl(\lim_{x\to\infty}e^{-ix\la}{\cal Z}_1(x,\la )\Bigr)
\Bigl(\lim_{x\to\infty}e^{ix\la}{\cal Z}_2(x,\la )\Bigr)^{-1},
\label{i27}\\ \noalign{\vskip6pt}
&&T_R(\la):=\Bigl(\lim_{x\to\infty}e^{ix\la}{\cal Z}_2(x,\la
)\Bigr)^{-1}.\label{i28}
\end{eqnarray}
Here the $p_1 \times p_1$ matrix ${\cal Y}_1$ and the $p_1
\times p_2$
matrix ${\cal Y}_2$ are upper and lower blocks of ${\cal Y}$,
respectively.
Analogously, the $p_1 \times p_2$ matrix ${\cal Z}_1$ and the
$p_2 \times
p_2$ matrix ${\cal Z}_2$ are upper and lower blocks of ${\cal Z}$,
respectively.
The functions $T_L$ and $T_R$ are called {\em the 
left} and {\em
the right transmission coefficients} and $R_L$ and
$R_R$ are
called {\em the  left} and {\em the right reflection
coefficients}, respectively.

Using (\ref{i24}) and  (\ref{i24''}) we get the result
\begin{Tm} \label{Tm2.2}  \cite{FKS} Let $\wt v\in$ GPE.
 Then the transmission and
reflection
coefficients are given by the formulas
\begin{eqnarray}
T_L(\la)&=&I_{p_1} +i\Phi_1(0)^*S_0^{-1}(\t - \lambda I_n 
)^{-1}
\Phi_1(0),
\label{i29}
\\
 R_L(\la)&=&i\Phi_2(0)^*S_0^{-1}(\t-\lambda I_n 
)^{-1} \Phi_1(0), \label{i30}
\\ 
 T_R(\la)&=& I_{p_2} + i\Phi_2(0)^*S_0^{-1}(\t-\lambda I_n 
)^{-1} (I_n- S_0 \om)\Phi_2(0) , \label{i31}
\\ \nonumber
 R_R(\la)&=& i\Phi_1(0)^*S_0^{-1}(\t-\lambda I_n 
)^{-1}(I_n-S_0 \om)\Phi_2(0)
\\ && 
+ i\Phi_1(0)^*(A^*-\lambda
I_n
)^{-1}\om\Phi_2(0), \label{i32}
\end{eqnarray}
where $ \theta: = A  - i\Phi_1(0)\Phi_1(0)^*S_0^{-1}$.
\end{Tm}
The inverse problem to recover $\wt v$ from $R_L$ is solved in Theorem 4.1
\cite{FKS}.
\begin{Tm} \label{InvRC}
Let ${\cal R}$ be a strictly proper rational $p_2 \times p_1$ matrix
function. Then   ${\cal R}$ is the left reflection coefficient of a
system $($\ref{i1}$)$ with $v \in {\mathrm{GPE}}$ if and only if
${\cal R}$  is  contractive on $\BR$. If ${\cal R}$ satisfies this condition,
then $\wt v$ can be uniquely recovered from ${\cal R}$ in two steps, that is,
steps (i) and (ii) below.
(i) First, take a minimal realization
\begin{equation}       \label{i33}
{\cal R}(\la)={\cal C}(\lambda I_n-{\cal A})^{-1}{\cal B}.
\end{equation}
Then there is a unique solution ${\cal X}$ of the Riccati equation
\begin{equation}\label{Ric}
i({\cal X}{\cal A}-{\cal A}^*{\cal X})={\cal C}^*{\cal C}+{\cal X}{\cal B}{\cal B}^*{\cal X}
\end{equation}
such that
\begin{equation}\label{i34}
\s({\cal A}+i{\cal B}{\cal B}^*{\cal X}) \subset \ov \BC_- , \quad {\cal X}={\cal X}^*.
\end{equation}
Moreover, $ \det {\cal X}\not=0$.  \\
(ii) Next, recover parameter matrices  by the equalities
\begin{equation}\label{i35}
A= {\cal A}+i  {\cal B} {\cal B}^* {\cal X}, \quad S_0= {\cal X}^{-1}, \quad \Phi_1(0)=  {\cal B}, \quad  
\Phi_2(0)= -i S_0
 {\cal C}^*,
\end{equation}
and recover $\wt v$ via $($\ref{i3}$)$, $($\ref{i4}$)$, and $($\ref{i22}$)$.
\end{Tm}
\subsubsection{Skew-self-adjoint Dirac system}
In this subsubsection we consider system
\begin{equation}       \label{i36}
\frac{d}{dx}\wt u(x, \la)=\big(i \la j+j \wt V(x)\big)\wt u(x,\la) \quad (0 \leq x <\infty),
\quad \wt V= \left[
\begin{array}{cc} 0 &
\wt v \\ \wt v^* & 0
\end{array}
\right],
\end{equation}
\begin{equation}       \label{i37}
\wt v(x)=2 \Phi_1(0)^*e^{ixA^*}S(x)^{-1}e^{ixA}\Phi_2(0), \quad
S(x)=I_n+\int_0^x \Pi(t)j\Pi(t)^*dt,
\end{equation}
where $j$ has the form $($\ref{r14}$)$, $\Pi$ is given by  $($\ref{i3}$)$, and the identity $A-A^*=i\Pi(0)\Pi(0)^*$
is fulfilled. That is, we  consider GBDT transformation of the trivial initial system, where $v \equiv 0$.
According to $($\ref{r11'}$)$ and to the second relation in  $($\ref{i37}$)$ 
the inequality $S(x)>0$ holds (see inequalities (1.6) and (1.7) in \cite{GKS2}).
Thus, the right-hand side of the first equality in $($\ref{i37}$)$ is well-defined.
The class of  potentials $\wt v$ of the form $($\ref{i37}$)$ is denoted by
PE$_2$. By Proposition 1.4 in \cite{GKS2} for  each $\wt v \in$ PE$_2$
there is some $M>0$ such that
\begin{equation}       \label{i38}
\sup_{x \in [0, \, \infty)}\|\wt v(x)\| \leq M.
\end{equation}

 It follows from Proposition \ref{PnGTskDS} that the fundamental solution
$\wt u(x,\la)$ of system  $($\ref{i36}$)$ (normalized by $\wt u(0,\la)=I_m$)  is given by formula $($\ref{i8}$)$,
where $w_A$ has the form $($\ref{se.63}$)$.

Analogously to Definition \ref{DnW1} one can define Weyl
functions of the skew-self-adjoint Dirac  system
\cite{ClGe, GKS2, SaA021, SaA14}.
\begin{Dn} \label{DnskW}
 Let  system  $($\ref{i36}$)$, where $\|\wt v(x)\|$  is bounded on all the finite intervals
$[0,\,l]$, be given.
  Then,
  a holomorphic $p \times p$ matrix function $\varphi$ such that
\begin{equation} \label{i39}
\int_0^\infty \left[ \begin{array}{lr}  \varphi (\la) ^* & I_p
\end{array} \right]
  \wt u(x, \la) ^*
\wt  u(x, \la) 
 \left[ \begin{array}{c}
\varphi (\la) \\ I_p  \end{array} \right] dx < \infty , \quad
\Im \la < -M_1
\end{equation}
is called a Weyl function  of  $($\ref{i36}$)$.
\end{Dn}
There is a unique Weyl function of system $($\ref{i36}$)$ with a bounded on $[0, \, \infty)$
potential  \cite{SaA021, SaA14}. If  $($\ref{i38}$)$ holds, we can put
in $($\ref{i39}$)$ $M_1=M$. For a  generalization of the Definition \ref{DnskW}  see \cite{ClGe, FKS3, SaAg}.

Direct and inverse problems for $\wt v \in$ PE$_2$ are solved explicitly
in Theorems 2.1 and 2.3 from \cite{GKS2}. It is done in a similar
to the self-adjoint case way.
\begin{Tm} \cite{GKS2} Let $v \in $
PE$_2$.
Then  system $($\ref{i36}$)$ has a unique Weyl function
$\vp$, which  satisfies $($\ref{i39}$)$ on the lower semi-plane
$ {\mathbb C}_{-}$, a finite number of poles excluded,
 and this function is given by the formula
\begin{equation}
\vp ( \lambda )=i \Phi_{1}(0)^{*}( \lambda I_{n} - \breve A )^{-1} \Phi_{2}(0), \quad \breve A:=A-i\Phi_{2}(0)\Phi_{2}(0)^*.
\end{equation}
 \end{Tm}
\begin{Tm}\label{InvSk} \cite{GKS2}
Let $ \vp $ be a strictly proper rational $p \times p$ matrix function.
Then $\vp$ is the Weyl function of the skew-self-adjoint Dirac system $($\ref{i39}$)$,
where $\wt v \in$ PE$_2$.
To recover $\wt v$ take a minimal realization
\begin{equation}       \label{i40}
\vp(\la)={\cal C}(\lambda I_n-{\cal A})^{-1}{\cal B}.
\end{equation}
There is a positive solution ${\cal X}>0$ of the  Riccati equation
\begin{equation}       \label{i41}
{\cal X}{\cal C}^*{\cal C}{\cal X}+i(   {\cal A}{\cal X} - {\cal X}{\cal A}^*)- {\cal B}{\cal B^*}=0.
\end{equation}
After putting
\begin{equation}       \label{i42}
\Phi_1(0)=i{\cal X}^{\frac{1}{2}}{\cal C}^*, \quad \Phi_2(0)={\cal X}^{-\frac{1}{2}}{\cal B}, \quad
 A={\cal X}^{-\frac{1}{2}}{\cal A}{\cal X}^{\frac{1}{2}}+i{\cal X}^{-\frac{1}{2}}{\cal B}{\cal B}^*
{\cal X}^{-\frac{1}{2}},
\end{equation}
the potential $\wt v$  is recovered by formula $($\ref{i37}$)$.
\end{Tm}
\subsection{$N$-wave equation}
Here we  consider GBDT transformation of the auxiliary system for
the $N$-wave equation
\begin{equation}       \label{i43}
\frac{d}{dx} \wt u(x,\la)=\big( i \la D-[D, \wt \xi(x) ]\big)u(x,\la)  \quad (0 \leq x <\infty), \quad
\wt \xi^*= \wt \xi .
\end{equation}
We consider the subcase of  $($\ref{se.66}$)$, where $B=I_m$, and assume that
\begin{equation}       \label{i44}
D={\mathrm{diag}}\, \{d_1,d_2, \ldots,d_m\}, \quad
d_1>d_2> \ldots>d_m>0.
\end{equation}
Starting from the initial system $($\ref{se.66}$)$ with $\xi \equiv 0$ and taking into account $B=I_m$, rewrite
 $($\ref{se.69}$)$ as:
\begin{equation}       \label{i45}
\wt \xi=\Pi^*S^{-1}\Pi,
\end{equation} 
where
\begin{equation}       \label{i46}
AS(0)-S(0)A^*=i\Pi(0)\Pi(0)^*, \quad      \Pi_x=-iA\Pi D, \quad S_x=\Pi D \Pi^*.
\end{equation} 
We require
\begin{equation}       \label{i46'}
\wt u(0,\la)=I_m, \quad S(0)>0.
\end{equation} 
As we have $u(x,\la)=\exp(ix\la D)$ for the fundamental solution of the
initial system, formula $($\ref{r19}$)$ takes the form
\begin{equation}       \label{i47}
\wt u(x,\la)=w_A(x, \la)e^{ix\la D}w_A(0, \la)^{-1},
\end{equation}
where $w_A$ has the form  $($\ref{se.63}$)$ and $AS-SA^*=i\Pi\Pi^*$.
It follows from $($\ref{se.63}$)$ and $AS-SA^*=i\Pi\Pi^*$ that
\begin{equation}       \label{i48}
w_A(x, \la)^*w_A(x, \la)=I_m-i(\la- \ov \la)\Pi(x)^*(A^*- \la I_n)^{-1}S(x)^{-1}(A- \la I_n)^{-1}\Pi(x).
\end{equation}
By the last relations in $($\ref{i46}$)$ and $($\ref{i46'}$)$ we get $S(x)>0$ for $x \geq 0$,
that is, $w_A(x,\la)$ is well-defined for all $x\geq 0$ and all $\la \not \in \s(A)$.
If $\la \in \BC_-$, the inequality $S(x)>0$ and formula $($\ref{i48}$)$ imply
\begin{eqnarray}       \label{i49}
&& w_A(x, \la)^*w_A(x, \la) \leq I_m , \\
&&     \label{i50}
i(\la- \ov \la)\Pi(x)^*(A^*- \la I_n)^{-1}S(x)^{-1}(A- \la I_n)^{-1}\Pi(x) \leq I_m .
\end{eqnarray}
\begin{Pn} \label{bnd} Let parameter matrices $A$, $S(0)$, and $\Pi(0)$ satisfy conditions
$AS(0)-S(0)A^*=i\Pi(0)\Pi(0)^*$ and $S(0)>0$. Then $\wt \xi$  determined by 
$($\ref{i45}$)$ and $($\ref{i46}$)$
is bounded on $[0, \, \infty)$.
\end{Pn}
\begin{proof}. It follows from $($\ref{i46}$)$ that $\Pi$ has the form
\begin{eqnarray}       \label{r21}
&&\Pi(x)=\Big[\exp \big(-i d_1x A\big) f_1 \,\,   \ldots \,\, \exp
\big(-i d_m xA\big) f_m\Big], \\ && \nonumber f_k\in \BC^m \quad (1 \leq k \leq m).
\end{eqnarray}
One can easily see
that ${\mathrm{span}}_{\la \in {\cal O}}(A- \la I_n)^{-1}f_k \supseteq f_k$ for any open domain ${\cal O}$.
Hence, by $($\ref{i50}$)$ and $($\ref{r21}$)$ we have
\begin{equation}       \label{r22}
\sup_{x \in \BR_+}\|S(x)^{-\frac{1}{2}}e^{-i d_kx A}f_k\|<\infty \quad (1 \leq k \leq m).
\end{equation}
Now, the boundedness of $\wt \xi$
given by $($\ref{i45}$)$  is immediate.
\end{proof}

The class of potentials $\wt \xi$ of the form $($\ref{i45}$)$, which are generated by parameter matrices
satisfying conditions of Proposition \ref{bnd}, is denoted by PE$_3$.

\begin{Dn} \cite{SaA5', SaAn} \label{NWeyl}
   A   Weyl function of system
(\ref{i43}) is an $m \times m$ matrix function $\varphi$, such
that for some $M>0$ it is analytic in a lower semi-plane $\Im \la <
-M$, and the inequalities
\begin{equation}  \label{i51}
      \sup\limits_{x \, \leq \, l, \ {\Im} \, \la \, <
\, -M} \,
          \bigl\|
                \wt u(x,\la)
                \varphi (\la)
                \exp \, (-i x \la D )
          \bigr\|
    < \infty
\end{equation}
hold for all $l < \infty$.
\end{Dn}
In our case $\wt \xi$ is bounded, and so  (see \cite{SaA5'}) the system  (\ref{i43})  has a unique
Weyl function $\wt \vp$ normalized by the condition
\begin{equation}
      \wt   \varphi_{ks}^{\,} (\la) \equiv 1  \, \mbox{
for   } \,  k = s, \quad
       \wt  \varphi_{ks}^{\,} (\la) \equiv 0 \, \mbox{
for   } \,  k > s.
    \label{r23}
\end{equation}
 This Weyl function satisfies for
some $\ve>0$ the inequality
\begin{equation}
      \int\limits_0^\infty
                   \Bigl(
                        \exp \, (i  x \overline  \la \,
D )
                   \Big)
                 \wt  \varphi (\la)^*
                   \wt u (x,\la)^*
                   \wt u(x,\la)
                   \wt \varphi (\la)
                  \exp \,  \Big(
                        x(- i \la D - \ve I_m^{\,} )
                        \Big)
      \, dx
    <  \infty.
    \label{r24}
\end{equation}
Here, we do not require  (\ref{r23}).
Our next theorem is immediate from $($\ref{i47}$)$ and $($\ref{i49}$)$.
\begin{Tm} \label{Ndir} A Weyl function
of system $($\ref{i43}$)$, where $\wt \xi \in {\mathrm{PE}}_3$,
is given by the formula
\begin{equation}       \label{i52}
\vp(\la)=w_A(0, \la)= I_m-i  \Pi(0)^*S(0)^{-1}(A-\la I_n)^{-1}\Pi(0).
\end{equation}
\end{Tm}
Notice that in view of   $($\ref{i48}$)$ and $($\ref{i49}$)$  the matrix function
$\vp(\la)$ given by $($\ref{i52}$)$ has  following properties:
\begin{equation}       \label{i53}
\vp(\la)\vp(\ov \la)^*= I_m, \quad \vp( \la)^*\vp( \la) \leq I_m \quad (\la \in \BC_-), \quad 
\lim_{\la \to \infty}\vp( \la)=I_m.
\end{equation}
In particular, this matrix function  admits a minimal realization
\begin{equation}       \label{i54}
\vp(\la)=I_m +{\cal C}(\lambda I_n-{\cal A})^{-1}{\cal B}.
\end{equation}
It follows from the second relation in $($\ref{i53}$)$ that 
$\s({\cal A}) \in \BC_+$. Hence, there is a unique and positive solution
$S_0>0$ of the identity
\begin{equation}       \label{i55}
{\cal A}S_0-S_0{\cal A}^*=i{\cal B}{\cal B}^*.
\end{equation}

Matrix functions admitting realization $($\ref{i53}$)$  satisfy conditions of Theorem 1.2 from \cite{SaA5'} and there is at most
one solution of the corresponding inverse problem.
This solution is given in Theorem 5.6 from \cite{SaAn}.
\begin{Tm} \label{Ninv} \cite{SaAn} Let a $p \times p$ rational matrix function
$\vp$ satisfy conditions $($\ref{i53}$)$. Then $\vp$ is a Weyl function
of the unique system 
$($\ref{i43}$)$, where $\wt \xi \in {\mathrm{PE}}_3$.
To recover $\wt \xi$ take a minimal realization $($\ref{i54}$)$ of $\vp$,
recover $S_0$ from $($\ref{i55}$)$, and put
\begin{equation}       \label{i56}
A={\cal A}, \quad S(0)=S_0, \quad \Pi(0)={\cal B}.
\end{equation}
Then $\xi$ is generated by the parameter matrices $A$, $S(0)$, and $\Pi(0)$
via formulas  $($\ref{i45}$)$ and $($\ref{i46}$)$.
\end{Tm}

Finally, let parameter matrices $A$, $S(0,0)>0$, and $\Pi(0,0)$
satisfy $($\ref{se.74}$)$, where $B=I_m$, and generate via $($\ref{i45}$)$
the solution $\wt \xi$ of the $N$-wave equation $($\ref{se.70}$)$.
The corresponding evolution of the Weyl function is given by the formula
\begin{equation}       \label{i57}
\vp(t,\la)=w_A(0,t, \la)= I_m-i  \Pi(0,t)^*S(0,t)^{-1}(A-\la I_n)^{-1}\Pi(0,t).
\end{equation}
Here $\Pi_t=-iA\Pi \wh D$, $S_t=\Pi \wh D\Pi^*$, and formula
$($\ref{i57}$)$ holds on the interval $t \in [0,\ve)$, where $S(0,t)>0$.

\subsection*{Acknowledgement}
The work was supported by the Austrian Science Fund (FWF) under
Grant  no. Y330.

\begin{appendix}

\section{Mathematical system theory}\label{mathsys}
\setcounter{equation}{0}

We present here some basic results and notions from
mathematical system theory of rational matrix functions
that are used in our review.
This material has its roots in  Kalman theory \cite{KFA}, and can be
found in various books (see, for instance, \cite{BGK1,CF}). See
also  interesting historical remarks in  \cite{VEK}.

The rational matrix functions appearing in the 
article are {\em proper}, that is, analytic at infinity. Such an $m_2
\times m_1$ matrix function $W$ can be represented in the form
\begin{equation}
\label{app.1} W(\lambda)={\cal D}+{\cal C}(\lambda I_n-{\cal A})^{-1}{\cal B},
\end{equation}
where ${\cal A}$ is a square matrix of some order $n$, the matrices ${\cal B}$
and ${\cal C}$ are of sizes $n \times m_1$ and $m_2 \times n$,
respectively, and ${\cal D}=W(\infty)$. The representation (\ref{app.1})
is called a {\it realization} of $W$, and the number $\mathrm{ord}({\cal A})$ (order of
the matrix ${\cal A}$) is called the {\it state space dimension} of the
realization.

 The realization
(\ref{app.1}) is said to be {\it minimal} if its state space
dimension $n$ is minimal among all possible realizations of $W$.
This minimal $n$  is called the {\it McMillan degree} of $W$. The
realization (\ref{app.1}) of $W$ is minimal if and only if
\begin{equation}
\label{app.2'} {\mathrm{span}}\bigcup_{k=0}^{n-1}\im
{\cal A}^k{\cal B}=\BC^n,\quad {\mathrm{span}} \bigcup_{k=0}^{n-1}\im  ({\cal A}^*)^k{\cal C}^*=\BC^n, \quad
n= \mathrm{ord}({\cal A}),
\end{equation}
where Im is image.
If for a pair of matrices ${\cal A}$, ${\cal B}$ the first equality in
(\ref{app.2'}) holds, then the pair ${\cal A}$, ${\cal B}$ is called {\it
controllable} or a {\it full range}. If the second equality in
(\ref{app.2'}) is fulfilled, then ${\cal C}$, ${\cal A}$ is said to be {\it
observable}. If a pair ${\cal A}$, ${\cal B}$
 is full range, and ${\cal K}$ is an $m_1 \times n$ matrix,  then
the pair ${\cal A}-{\cal B}{\cal K}$, ${\cal B}$ is also full range.

Minimal realizations are unique up to a basis transformation, that
is, if (\ref{app.1}) is a minimal realization of $W$ and if
$W(\lambda)={\cal D}+\widetilde {\cal C}(\lambda I_n -\widetilde
{\cal A})^{-1}\widetilde {\cal B}$ is a second minimal realization of $W$, then
there exists an invertible matrix ${\cal S}$ such that
\begin{equation}
\label{app.3} \widetilde{\cal A}={\cal S}{\cal A}{\cal S}^{-1},\quad \widetilde {\cal B}={\cal S}{\cal B},\quad
\widetilde {\cal C}={\cal C}{\cal S}^{-1}.
\end{equation}
In this case, (\ref{app.3})  is called a {\it similarity
transformation}.

Finally, if $W$ is a square matrix and ${\cal D}=I_{m_1}$ ($m_1=m_2$), then $W^{-1}$
admits representation
\begin{equation} \label{app4}
W(\lambda)^{-1}=I_{m_1}-{\cal C}(\lambda I_n-{\cal A}^\times)^{-1}{\cal B},\quad
{\cal A}^\times={\cal A}-{\cal B}{\cal C}.
\end{equation}
\end{appendix}

\begin{flushright} \it
A.L. Sakhnovich, \\  Fakult\"at f\"ur Mathematik,
Universit\"at Wien,
\\
Nordbergstrasse 15, A-1090 Wien, Austria \\ 
e-mail: {\tt al$_-$sakhnov@yahoo.com }
\end{flushright}

\end{document}